\documentclass[a4paper,reqno]{amsart}

\usepackage[utf8x]{inputenc}
\usepackage[T1]{fontenc}
\usepackage{times,eulervm}
\usepackage{enumerate}
\usepackage{booktabs}
\usepackage{hyphenat}
\usepackage{amssymb}
\usepackage{siunitx}
\usepackage{bm}
\usepackage{hyperref}
\usepackage[capitalize]{cleveref}% After hyperref!
\usepackage{subcaption}
\usepackage[super]{nth}
\usepackage{algorithm}
\usepackage{algorithmic}
\usepackage{placeins}

\usepackage{tikz}
\usepackage{stanli}
\usetikzlibrary{decorations.pathreplacing}
\usepackage{blkarray}
\usepackage{pgfplots}
\pgfplotsset{compat=1.18}
\usetikzlibrary{pgfplots.statistics, pgfplots.colorbrewer}
\usepackage{pgfplotstable}

\usepackage[ngerman,american]{babel}% Für temporäres ...
\usepackage{ao-abk}% ... Einbinden von entwurf-text.tex.

\usepackage{ao-abbr}
\usepackage{ao-ximg}% Includes graphicx, luh-colors.
\usepackage{ao-matrix}
\usepackage[sepmid]{ao-math-opt}% Optimization; includes ao-math.std.
\usepackage{ao-math-fields}
\usepackage{ao-math-symbols}
\usepackage{ao-code-names}

\newcommand\bB{\bm B}
\newcommand\bK{\bm K}
\newcommand\bM{\bm M}
\newcommand\bY{\bm Y}
\newcommand\bS{\bm S}
\newcommand\bb{\bm b}
\newcommand\bu{\bm u}
\newcommand\bw{\bm w}
\newcommand\bx{\bm x}
\newcommand\bPhi{\bm\Phi}
\newcommand\ddf[2]{\Delta\ifgiven{#1}{(#1; #2)}{}}% Damage distr. fct.
\newcommand\cpdf[2]{F\ifgiven{#1}{(#1; #2)}{}}% Cumulative. distr. fct.
\newcommand\sbM[2][]{_{\itxt{M}#2\optindex{#1}}}% Subindex M#2[,#1].
\newcommand\sbS[2][]{_{\itxt{S}#2\optindex{#1}}}% Subindex S#2[,#1].

\newcommand\setM{\mathcal M}
\newcommand\setN{\mathcal N}
\renewcommand\bigO{\mathcal O}
\newcommand\Eff{\textup{Eff}}
\newcommand\ND{\textup{ND}}
\newcommand\GYRM{\textup{GraefYounesReductionMethod}}
\newcommand\uHoF{\textup{updateHallOfFame}}

\begin{document}

\title[Damage Location by Multi-Objective Pattern Search]
{Damage Location in Mechanical Structures\\
  by Multi-Objective Pattern Search}

\author[Günther]{Christian Günther}
\address{Christian Günther\\
  Leibniz Universität Hannover\\Institute of Applied Mathematics\\
  Wel\-fen\-gar\-ten 1\\30167 Hannover\\Germany}
\email{c.guenther@ifam.uni-hannover.de}
\urladdr{ifam.uni-hannover.de/de/cguenther}

\author[Hofmeister]{Benedikt Hofmeister}
\address{Benedikt Hofmeister\\
	Leibniz Universität Hannover\\Institute of Structural Analysis\\
	Ap\-pel\-straße 9A\\30167\\Hannover\\Germany}
\email{b.hofmeister@isd.uni-hannover.de}
\urladdr{https://www.isd.uni-hannover.de/de/hofmeister}

\author[Hübler]{Clemens Hübler}
\address{Clemens Hübler\\
  TU Darmstadt\\Institute of Structural Mechanics and Design\\Fran\-zis\-ka-Braun-Straße 3\\
  64287\\Darmstadt\\Germany}
\email{huebler@ismd.tu-darmstadt.de}
\urladdr{https://www.ismd.tu-darmstadt.de/das\_institut\_ismd/mitarbeiter\_innen\_ismd/team\_ismd\_details\_130944.de.jsp}

\author[Jonscher]{Clemens Jonscher}
\address{Clemens Jonscher\\
	Leibniz Universität Hannover\\Institute of Structural Analysis\\
	Ap\-pel\-straße 9A\\30167\\Hannover\\Germany}
\email{c.jonscher@isd.uni-hannover.de}
\urladdr{https://www.isd.uni-hannover.de/de/jonscher}

\author[Ragnitz]{Jasper Ragnitz}
\address{Jasper Ragnitz\\
	Leibniz Universität Hannover\\Institute of Structural Analysis\\
	Ap\-pel\-straße 9A\\30167\\Hannover\\Germany}
\email{j.ragnitz@isd.uni-hannover.de }
\urladdr{https://www.isd.uni-hannover.de/de/institut/personenverzeichnis/jasper-ragnitz-msc}

\author[Schubert]{Jenny Schubert}
\address{Jenny Schubert\\
  Leibniz Universität Hannover\\Institute of Applied Mathematics\\
  Wel\-fen\-gar\-ten 1\\30167 Hannover\\Germany}
\email{schubert@ifam.uni-hannover.de}
\urladdr{ifam.uni-hannover.de/schubert}

\author[Steinbach]{Marc C. Steinbach}
\address{Marc C. Steinbach\\
  Leibniz Universität Hannover\\Institute of Applied Mathematics\\
  Wel\-fen\-gar\-ten 1\\30167 Hannover\\Germany}
\email{mcs@ifam.uni-hannover.de}
\urladdr{ifam.uni-hannover.de/mcs}

\begin{abstract}
  % Introduction. What’s the topic?
  We propose a multi-objective global pattern search algorithm
  for the task of locating and quantifying damage
  in flexible mechanical structures.
  % State the problem you tackle.
  This is achieved by identifying eigenfrequencies and eigenmodes
  from measurements and matching them
  against the results of a finite element simulation model,
  which leads to a nonsmooth nonlinear
  bi-objective parameter estimation problem.
  % Summarize why nobody else has adequately answered the research question yet.
  A derivative-free optimization algorithm is required
  since the problem is nonsmooth
  and also because complex mechanical simulation models
  are often solved using commercial black-box software.
  Moreover, the entire set of non-dominated solutions
  is of interest to practitioners.
  Most solution approaches published to date
  are based on meta-heuristics such as genetic algorithms.
  % Explain how you tackled the research question.
  The proposed multi-objective pattern-search algorithm
  provides a mathematically well-founded alternative.
  It features a novel sorting procedure
  that reduces the complexity in our context.
  % How did you go about doing the research that follows from your big idea.
  Test runs on two experimental structures
  with multiple damage scenarios
  are used to validate the approach.
  % What’s the key impact of your research?
  The results demonstrate that the proposed algorithm
  yields accurate damage locations
  and requires moderate computational resources.
  From the engineer's perspective
  it represents a promising tool for structural health monitoring.
\end{abstract}


% MCS: hier ist die Auswahl schwierig ...
\subjclass[2020]{%
  %35R30, % Inverse problems for PDEs
  %65K05, % Numerical mathematical programming methods
  %74G75, % Inverse problems in equilibrium solid mechanics
  70J10, % Modal analysis in linear vibration theory
  %74K10, % Rods (beams, columns, shafts, arches, rings, etc.)
  % 74Pxx: Optimization problems in solid mechanics
  %74P99, % None of the above, but in this section
  % 74Rxx: Fracture and damage
  74R99, % None of the above, but in this section
  90C26, % Nonconvex programming, global optimization
  90C29, % Multi-objective and goal programming
  %90C30, % Nonlinear programming
  %90C56, % Derivative-free methods and methods using generalized derivatives
  90C59%, % Approximation methods and heuristics in mathematical programming
  %90C90, % Applications of mathematical programming
  %92F05, % Other natural sciences (mathematical treatment)
}

\date\today

\maketitle

\section{Introduction}
A common task in structural engineering
is the updating of finite element (FE) models
wherein measurements of the mechanical structure
are used to estimate parameters of the simulation model
\cite{friswell_finite_1995, mottershead1993}.
This is necessary when physical properties and hence model parameters
may change over time due to wearout, damage, or other reasons.
The present work, in particular, addresses the problem
of locating and also quantifying damage in terms of stiffness deviations
based on detected changes in eigenfrequencies and mode shapes,
which is known as \emph{vibration-based model updating}
in structural health monitoring \cite{fan2011}.

Mathematically this task leads to a nonsmooth and nonlinear
bi-objective least-squares parameter estimation problem.
Changes of eigenmodes and eigenfrequencies
are treated as separate error measures
to obtain more error-tolerant results.
These error measures depend on certain damage parameters
that influence the stiffness matrix
and possibly the inertia matrix
of the mechanical structure.
As the eigenmodes and eigenfrequencies
are solutions of an eigenvalue problem for the FE simulation model,
their dependence on the parameters is not only nonsmooth and nonlinear
but also expensive to evaluate.
Moreover, the entire set of non-dominated solutions
(the Pareto front or Pareto frontier)
is of interest for assessing the outcome
of the model updating problem
in practical applications.
We propose a multi-objective global pattern search algorithm
with novel algorithmic elements
to compute large numbers of non-dominated points
close to the Pareto front.

Finite element model updating was first proposed in the 1990s,
notable early contributors are Friswell and Mottershead
\cite{friswell_finite_1995} as well as Link \cite{Link1999}.
Commonly, metaheuristic optimization algorithms are used for model updating
\cite{levin_dynamic_1998, jahjouh_modified_2016, jahjouh2016}.
Other approaches make use of local optimization algorithms
\cite{petersen_2017} or a combination of a global and a local optimizer
\cite{schroeder2018, begambre_hybrid_2009}.
Moreover, previous approaches to model updating
are frequently based on convex combinations
of the mode shape and eigenfrequency error metrics \cite{Simoen2015}.
However, since there is no natural ``best'' weighting
for a given application problem
and an appropriate weighting is hard to guess
even for experienced practitioners,
a priori choices like equal weights usually yield inadequate results
\cite{schroeder2018}.
The proper mathematical tool to deal with this situation
is multi-objective optimization \cite{Naranjo2020,custodio2012}.

Various algorithmic approaches have been developed
for solving specific classes of multi-objective optimization problems.
These include the generalization of well-known derivative-based
single-objective optimization methods, in particular line search algorithms
such as the steepest descent method and the Newton algorithm;
see, \eg, \cite{FliegeEtAl2000,DrummondSvaiter2005,FliegeEtAl2009}.
Although derivative-based methods are commonly preferred
over derivative-free methods because of their superior performance,
this view may change if the goal is to compute an approximation to the
entire Pareto front of a given multi-objective optimization problem,
as recognized in \cite{CustodioEtAl2022}.
It is well-known that methods based on solving
parametic families of scalar optimization problems
(so-called scalarization approaches; e.g. \cite{Ehrgott2005,Jahn11})
are useful for the approximation of the Pareto front.
Moreover, in recent years also branch and bound methods
have been proposed to obtain approximations (coverages)
of the set of Pareto efficient solutions and of the Pareto front;
\eg, \cite{EichfelderNiebling2019,EichfelderEtAl2021a,EichfelderEtAl2021b}.

Similarly, derivative-free methods are often the preferred choice
if a multi-objective optimization problem involves
highly expensive objective functions.
In \cite{ThomannEichfelder2020} a trust region type algorithm is presented
for multi-objective optimization problems in which one objective function
is significantly more expensive then the remaining ones,
though under differentiability assumptions on all involved functions.
In the model updating application considered here,
all objective functions are expensive and nondifferentiable;
hence we cannot emply derivative-based algorithms
to approximate the desired Pareto fronts.
This situation is common in engineering applications
where complex numerical models are considered.

Derivative-free methods are generally classified as
deterministic, non-determinisic or hybrid.
On the same input, deterministic methods always produce the same output,
non-deterministic ones may produce different output in different runs,
and hybrid methods combine deterministic and non-deterministic features.
Deterministic derivative-free methods for multi-objective optimization
(in particular for approximating the entire Pareto front)
have attracted much less attention than non-deterministic ones.
Some well-known methods of this type are proposed in
\cite{custodio2011,custodio2012,evtushenko2014,alotto2015};
some hybrid methods are proposed in
\cite{Jahn2006,custodio2011,alotto2015}.
In particular, pattern search (or direct search) approaches
as extensions of the classical pattern search
play a prominent role in this field; \eg,
\cite{Hooke1961,torczon1997,Dolan2003,abramson2005}.

Non-deterministic derivative-free methods are often used in engineering
because they are comparatively easy to implement and to use.
Methods of this type follow a random-number-based approach
and are used in evolutionary (population-based) multi-objective optimization;
see \cite{custodio2012} for a survey on this topic.
A particularly prominent example
is the multi-objective genetic algorithm NSGA-II \cite{deb2000}.
A more recent branch of meta-heuristic multi-objective optimizers
includes particle swarm approaches, often referred to as MOPSO;
see \cite{Coello2004} for a survey.
We finally notice that there is a very large body
of literature on related non-deterministic methods which,
however, are not of primary interest in our context
as we prefer a mathematically better founded approach.

In this work we propose a deterministic derivative-free algorithm,
based on a pattern search approach,
to compute large numbers of non-dominated points close to the Pareto front
of the given multi-objective optimization problem.
To this end, a novel non-dominated sorting procedure
for computing first level and also higher level Pareto fronts
(similar to NSGA-II in \cite{deb2000})
will play a key role in our algorithm.

The remainder of the paper is structured as follows.
In \cref{sec:application} we provide some background information
on the application problem.
The precise mathematical formulations of two typical FE simulation models,
of our damage model and of the resulting bi-objective optimization problem
are presented in \cref{sec:model}.
In \cref{sec:algorithm} we describe and discuss the details
of our proposed multi-objective global pattern search algorithm.
In \cref{sec:results} the algorithm is then tested
on two mechanical structures with artificial damage
to evaluate and demonstrate its effectiveness.
We conclude the article with some final remarks and an outlook
in \cref{sec:conclusion}.

%%% Local Variables:
%%% mode: latex
%%% TeX-master: t
%%% End:

\section{Application Problem}
\label{sec:application}
In this section we give some additional information
on the model updating application.
Mathematical details are provided in the following section.
While the basic idea of modifying an FE simulation model
to match observed changes of the structural behavior
always leads to a parameter estimation problem,
there is a wide range of situations where model updating is relevant
and a multitude of possibilities how to approach this task.
The range of situations regards the types and sizes of structures,
the types of possible property changes, failures, or damage,
the type and quality of the measured data,
and the type and amount of derived results required by the engineer
to determine the new model parameters.
The possibilities of approaching the task may differ
in the parameterization of the FE model,
in the choice of objective functions and possibly constraints
(thus in the mathematical properties of the optimization problem),
and in the choice of the solution algorithm.

As already indicated, we address the problem of locating and quantifying
damage in flexible structures.
More precisely, we restrict ourselves to structures
that can be modeled as ensembles of one-dimensional flexible beams,
and we consider two beam types of different complexity,
the Euler-Bernoulli beam and the Timoshenko beam.
In reality the structures of interest may be very large,
such as rotor blades of wind turbines.
The practicability and cost of obtaining useful measurements determines
which structural properties can be monitored in real life.
Since acceleration sensors are small and cheap,
and since they can continuously deliver data at high sampling rates,
we focus on time series of acceleration measurements
as data basis for the desired damage location.
There are standard methods to determine eigenfrequencies
and associated mode shapes of the structure from this data.
In principle, we are thus able to locate any kind of damage
that changes the stiffness or the inertia of the structure,
by choosing changes of the eigenfrequencies and of the mode shapes
as two distinct error measures in the parameter estimation problem.
This gives the mentioned nonlinear and nonsmooth bi-objective
parameter estimation problem, where the choice of constraints
depends on the chosen model parameterization
and also on the specific application problem.
The proposed multi-objective global pattern-search algorithm
is then a suitable choice for solving the problem
because the entire Pareto front is of interest
for determining the updated model parameters.

A central mathematical aspect deserves further discussion:
the parameterization of the FE model.
The most straightforward option is to assign one parameter
to each individual finite element \cite{levin_dynamic_1998}.
As FE models often have hundreds or thousands of elements,
this approach leads to an excessive number of parameters
and typically to an ill-posed optimization problem \cite{Simoen2015}.
To reduce the number of parameters,
a strategy commonly used is the assignment of one parameter
to a group of elements supposedly having similar mechanical properties,
which are then referred to as substructures or super-elements
\cite{kim_improved_2004}.
The more super-elements are employed,
the higher the number of parameters becomes.
A different approach, which is motivated by a variation
of the properties along the structure, uses a damage distribution.
This also leads to a formulation with comparatively few parameters.
The distribution functions can take various forms,
for example quadratic polynomials \cite{teughels_global_2003}.
In the present work we employ Gaussian damage distribution functions
as suggested in \cite{Bruns2019_1}
to parameterize the structural stiffness of the FE model.

\section{Mathematical Model}
\label{sec:model}
\subsection{Structural Model}
\label{sec:model:struct}

\begin{figure}
  \centering
  \begin{tikzpicture}
  %the points
  \point{begin}{0}{0};
  \point{end0}{0}{0.4};
  \point{end1}{3}{0.4};
  \point{end2}{6}{0.4};
  \point{end3}{8}{0.4};
  \point{end4}{11}{0.4};
  \point{end0n}{0}{0.65};
  \point{end1n}{3}{0.65};
  \point{end2n}{6}{0.65};
  \point{end3n}{8}{0.65};
  \point{end4n}{11}{0.65};
  \node at (0,-1.4) {$s_0$};
  \node at (3,-1.4) {$s_1$};
  \node at (6,-1.4) {$s_2$};
  \node at (8,-1.4) {$s_{n-1}$};
  \node at (11,-1.4) {$s_n$};
  %support
  \support{3}{begin}[-90];
  \load{2}{end0}[180][-180];
  \notation{1}{end0n}{$\psi_{11}$}[above right]
  \load{2}{end1}[180][-180];
  \notation{1}{end1n}{$\psi_{12}\enspace\psi_{21}$}[above]
  \load{2}{end2}[180][-180];
  \notation{1}{end2n}{$\psi_{22}\enspace\psi_{31}$}[above]
  \load{2}{end3}[180][-180];
  \notation{1}{end3n}{\kern-1.2em$\psi_{n-1,2}\enspace\psi_{n1}$}[above]
  \load{2}{end4}[180][-180];
  \notation{1}{end4n}{$\psi_{n2}$}[above left]
  \dimensioning{1}{end0}{end1}{-1}[$l_1$];
  \dimensioning{1}{end1}{end2}{-1}[$l_2$];
  \dimensioning{1}{end3}{end4}{-1}[$l_n$];
  %beams
  \node[rectangle,draw, anchor = west, line width=0.4mm, minimum width=3cm, minimum height=0.7cm] (r) at (0,0) {$E_1, I_1, A_1, \rho_1$};
  \node[rectangle,draw, anchor = west, line width=0.4mm, minimum width=3cm, minimum height=0.7cm] (r) at (3,0) {$E_2, I_2, A_2, \rho_2$};
  \node[rectangle,draw, anchor = west, line width=0.4mm, minimum width=2cm, minimum height=0.7cm] (r) at (6,0) {};
  \node[rectangle,draw, anchor = west, line width=0.4mm, minimum width=1cm, minimum height=1cm, color = white, fill = white] (r) at (6.5,0) {};
  \node[anchor = center] (r) at (7,0) {\dots};
  \node[rectangle,draw, anchor = west, line width=0.4mm, minimum width=3cm, minimum height=0.7cm] (r) at (8,0) {$E_n, I_n, A_n, \rho_n$};
  %DoFs
  \draw[->,line width=0.03] (0.015,0) -- (0.015,1.5) node[right] {$u_{11}$};
  \draw[->,line width=0.03] (3.015,0) -- (3.015,1.5) node {$u_{12}\enspace u_{21}$};
  \draw[->,line width=0.03] (6.015,0) -- (6.015,1.5) node {$u_{22}\enspace u_{31}$};
  \draw[->,line width=0.03] (8.015,0) -- (8.015,1.5) node {\kern-1.2em$u_{n-1,2}\enspace u_{n1}$};
  \draw[->,line width=0.03] (11.015,0) -- (11.015,1.5) node[left] {$u_{n2}$};
\end{tikzpicture}
  \caption{Scheme of a cantilever beam with $n$ finite elements.}
  \label{fig:beam}
\end{figure}
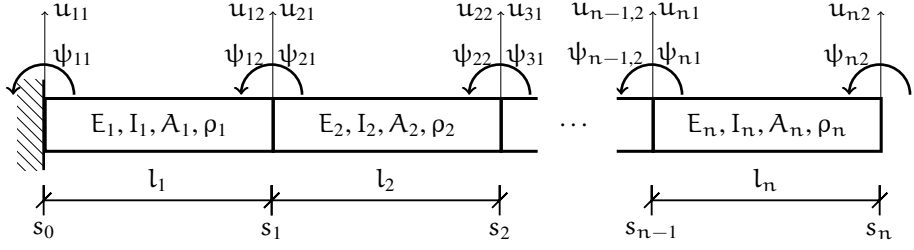

We focus on slender beam structures,
considered as essentially one-dimensional mechanical objects.
The mathematical models follow a standard finite element approach.
A schematic view of a discretized beam is given in \cref{fig:beam}.

For the applications considered later, we need two types of models:
the Euler-Bernoulli beam and the Timoshenko beam.
We start with the simpler Euler-Bernoulli beam model
\cite[Chapter~5.4]{friswell2010}
wherein shear effects and rotational inertia are ignored.
The beam of length $L$ is discretized into a sequence of $n$ finite elements
(not necessarily of equal length)
with node positions $s_i \in [0, L]$, $i \in \setN = \set{0, \dots, n}$.
The degrees of freedom (DoF) of each beam element
$e \in \setE = \set{1, \dots, n}$ are located
at the two end nodes $e - 1, e$, locally indexed as $e1, e2$.
Each node has two DoF: a lateral displacement $u_{ei}\in\R$
and a rotation angle $\psi_{ei} \in \R$, $i \in \set{1, 2}$.
An approximate displacement of a point $\xi \in [0, l_e]$
on beam element $e$ at time $t$ is calculated with four shape
functions $N_{ei}\: [0, l_e] \to \R$ as well as the four DoF:
\begin{equation}
  \label{eq:displacementEBB}
  u_e(\xi,t)
  =
  \mat[4]{N_{e1}(\xi) & N_{e2}(\xi) & N_{e3}(\xi) & N_{e4}(\xi)} q_e(t),
  \quad
  q_e(t)
  =
  \col(c){u_{e1}(t) \\ \psi_{e1}(t) \\ u_{e2}(t) \\ \psi_{e2}(t)},
\end{equation}
where the shape functions are defined as
\begin{equation}
\label{eq:shape_functions}
\begin{aligned}
  N_{e1}(\xi)
  &= \left( 1 - 3 \frac{\xi^2}{l_e^2} + 2 \frac{\xi^3}{l_e^3} \right),
  \quad &
  N_{e2}(\xi)
  &= l_e \left(
    \frac{\xi}{l_e} - 2 \frac{\xi^2}{l_e^2} + \frac{\xi^3}{l_e^3} \right),
  \\
  N_{e3}(\xi)
  &= \left( 3 \frac{\xi^2}{l_e^2} - 2 \frac{\xi^3}{l_e^3} \right),
  &
  N_{e4}(\xi)
  &= l_e \left( -\frac{\xi^2}{l_e^2} + \frac{\xi^3}{l_e^3} \right).
\end{aligned}
\end{equation}
Here $l_e = s_e - s_{e-1}$ is the length of beam element $e$,
and \eqref{eq:displacementEBB} is a cubic polynomial in $\xi$.

The stiffness matrix $\bK_e \in \R^{4 \x 4}$
and the mass matrix $\bM_e \in \R^{4 \x 4}$
of beam element $e$ are derived from (approximations of)
the strain energy $U_e$ and the kinetic energy $T_e$.
Under the assumptions of a linear elastic material
and a constant area moment of inertia $I_e$ (also called second moment of area),
the strain energy is approximated as
\begin{equation*}
  U_e(t)
  =
  \frac12 E_e I_e \int_0^{l_e} u_e''(\xi, t)^2 d\xi,
\end{equation*}
where $E_e$ is Young's modulus and $u_e'$ denotes
the partial derivative with respect to $\xi$.
The area moment of inertia is
\begin{equation*}
  I_e = \dfrac{w_e t_e^3}{12},
\end{equation*}
where $w_e$ is the width of the rectangular beam cross section
and $t_e$ is its thickness \cite[Appendix~A.1]{timoshenko1956}.
Together with \eqref{eq:displacementEBB} we obtain
\begin{equation*}
  U_e(t) = \frac12 q_e(t)^T \bK_e q_e(t),
\end{equation*}
with the elements of the stiffness matrix being
\begin{equation*}
  k_{e,ij} = E_e I_e \int_0^{l_e} N_{ei}''(\xi) N_{ej}''(\xi) \, d\xi.
\end{equation*}
The stiffness matrix of beam element $e$ then evaluates to
\begin{equation}
  \label{eq:ElementStiffness}
  \bK_e
  =
  \frac{E_e I_e}{l_e^3}
  \mat[4]{
    12    &  6 l_e    & -12   &  6 l_e \\
    6 l_e &  4 l_e ^2 & -6 l_e &  2 l_e ^2 \\
    -12   & -6 l_e    & 12     & -6 l_e  \\
    6 l_e &  2 l_e ^2 & -6 l_e &  4 l_e ^2
  }.
\end{equation}
With cross-sectional area $A_e = w_e t_e$,
the kinetic energy of beam element $e$ is approximated as
\begin{equation*}
  T_e(t) = \frac12 \rho_e A_e \int_0^{l_e} \dot u_e(\xi, t)^2 \, d\xi,
\end{equation*}
where $\rho_e$ denotes the (constant) density of the material
and $\dot u_e$ denotes the partial derivative with respect to $t$.
With \eqref{eq:displacementEBB} we obtain
\begin{equation*}
  T_e(t) = \frac12 \dot q_e(t)^T \bM_e \dot q_e(t),
\end{equation*}
where the elements of the mass matrix are
\begin{equation*}
  m_{e,ij} = \rho_e A_e \int_0^{l_e} N_{ei}(\xi) N_{ej}(\xi) \, d\xi.
\end{equation*}
The mass matrix of beam element $e$ then evaluates to
\begin{equation}
  \label{eq:ElementMass}
  \bM_e
  =
  \frac{l_e\rho_e A_e}{420}
  \mat[4]{
    156     & 22 l_e   &  54    & -13 l_e \\
    22 l_e  &  4 l_e^2 &  13 l_e &  -3 l_e^2 \\
    54      & 13 l_e   & 156     & -22 l_e \\
    -13 l_e & -3 l_e^2 & -22 l_e &   4 l_e^2
  }.
\end{equation}
The stiffness matrix and mass matrix are both symmetric positive definite.

Since each inner node $i \in \set{1, \dots, n - 1}$
belongs to two beam elements, $e = i$ and $e = i + 1$,
both elements contribute to the mass matrix of the entire beam.
The latter has dimension $2 (n + 1)$ and is
symmetric positive definite and block-tridiagonal,
\begin{equation*}
  \bM
  =
  \mat[6]{
    M_1^{11} & M_1^{12} \\
    M_1^{21} & \_M_2   & M_2^{12} \\
            & M_2^{21} & \_M_3   & M_3^{12} \\
            &         & M_3^{21} & \ddots & \ddots \\
            &         &         & \ddots & \_M_n   & M_n^{12} \\
            &         &         &        & M_n^{21} & M_n^{22}
          },
\end{equation*}
where $\_M_i = M_{i-1}^{22} + M_i^{11}$ and
\begin{equation*}
  M_e^{ij}
  =
  \mat[2]{
    m_{e,2i-1,2j-1} & m_{e,2i-1,2j} \\
    m_{e,2i,2j-1}   & m_{e,2i,2j}},
  \quad i,j \in \set{1, 2}.
\end{equation*}
The full stiffness matrix $\bK$ has precisely the same structure.
In the common special case where the beam is fixed at node $0$,
the translational and rotational displacements are
\begin{equation*}
  u_{11}(t) = 0, \quad \psi_{11}(t) = 0.
\end{equation*}
Thus the first two DoFs do not appear in the finite element model,
and we delete the first two rows and columns to obtain
$\bM, \bK \in \R^{2 n \x 2 n}$ (without changing notation).

When considering a damaged beam,
we assume that the mass matrix $\bM$ does not change
whereas the stiffness matrix $\bK(\bx)$ depends on
the parameter vector $\bx$ of the damage model.
This dependence is assumed to be caused by a dependence
of each Young's modulus on $\bx$,
with details given in \cref{sec:model:damage},
\begin{equation}
  \label{eq:ElementStiffness_Damaged}
  \bK_e(x)
  =
  \frac{E_e(x) I_e}{l_e^3}
  \mat[4]{
    12    &  6 l_e    & -12   &  6 l_e \\
    6 l_e &  4 l_e ^2 & -6 l_e &  2 l_e ^2 \\
    -12   & -6 l_e    & 12     & -6 l_e  \\
    6 l_e &  2 l_e ^2 & -6 l_e &  4 l_e ^2
  }.
\end{equation}

To account for shear effects and rotational inertia,
we now introduce the extended model of a Timoshenko beam
\cite[Chapter~5.4]{friswell2010}.
There are approaches that increase the DoF per node to account for shear.
However, we follow a reduced integration approach
to retain only two DoF per node:
the lateral displacement $u_e$
and the angle of the beam cross section $\psi_e$.
The angle of the beam cross section is given as
\begin{equation*}
  \psi_e(\xi,t) = u_e'(\xi,t) + \beta_e(\xi,t).
\end{equation*}
Herein $\beta_e$ denotes the shear angle,
\begin{equation}
  \label{eq:shear-angle}
  \beta_e(\xi,t) = \frac{\Phi_el_e^2}{12} u_e'''(\xi,t),
  \qquad
  \Phi_e = \frac{12E_eI_e}{\kappa_eG_eA_el_e^2},
\end{equation}
where $\kappa_e$ is the shear constant
(which depends on the shape of the cross section)
and $G_e$ is the shear modulus.
Similar to the Euler-Bernoulli beam, we obtain
\begin{equation}
  \label{eq:displacementEBB_timo}
  u_e(\xi,t)
  =
  \mat[4]{N_{e1}(\xi) & N_{e2}(\xi) & N_{e3}(\xi) & N_{e4}(\xi)} q_e(t),
\end{equation}
where $q_e$ is defined as in \eqref{eq:displacementEBB}
and the shape functions are now defined as
\begin{equation}
\label{eq:shape_functions_timo}
\begin{aligned}
  N_{e1}(\xi)
  &= \frac{1}{1+\Phi_e} \left( 1 + \Phi_e - \Phi_e \frac{\xi}{l_e}
    - 3 \frac{\xi^2}{l_e^2} + 2\frac{\xi^3}{l_e^3} \right), \\
  N_{e2}(\xi)
  &= \frac{l_e}{1+\Phi_e} \left( \frac{2+\Phi_e}{2} \frac{\xi}{l_e}
    - \frac{4+\Phi_e}{2} \frac{\xi^2}{l_e^2} + \frac{\xi^3}{l_e^3} \right), \\
  N_{e3}(\xi)
  &= \frac{1}{1+\Phi_e} \left( \Phi_e \frac{\xi}{l_e}
    + 3\frac{\xi^2}{l_e^2} - 2 \frac{\xi^3}{l_e^3} \right), \\
  N_{e4}(\xi)
  &= \frac{l_e}{1+\Phi_e} \left( -\frac{\Phi_e}{2} \frac{\xi}{l_e}
    - \frac{2-\Phi_e}{2} \frac{\xi^2}{l_e^2} + \frac{\xi^3}{l_e^3} \right).
\end{aligned}
\end{equation}
If $\Phi_e = 0$, shear effects are neglected
and the shape functions \eqref{eq:shape_functions_timo}
are identical to \eqref{eq:shape_functions}.
The strain energy for the element, including shear, is
\begin{equation*}
  U_e
  =
  \frac12 E_e \int_0^{l_e} I_e(\xi) \psi_e'(\xi,t)^2 \, d\xi
  +
  \frac12 \kappa_e G_e \int_0^{l_e} A_e(\xi)\beta_e(\xi,t)^2 \, d\xi.
\end{equation*}
For a constant element cross section we have
\begin{equation*}
  \psi_e'(\xi,t) = u_e''(\xi,t) + \beta_e'(\xi,t) = u_e''(\xi,t),
\end{equation*}
because $\beta_e$ is constant along the element.
Thus, the approximate strain energy is
\begin{equation*}
  U_e(t) = \frac12 q_e(t)^T \bK_e q_e(t),
\end{equation*}
with the elements of the stiffness matrix being
\begin{equation*}
  k_{e,ij}
  =
  E_e I_e \int_0^{l_e} N_{ei}''(\xi) N_{ej}''(\xi) \, d\xi
  +
  \frac{E_e I_e \Phi_e l_e^2}{12}
  \int_0^{l_e} N_{ei}'''(\xi) N_{ej}'''(\xi) \, d\xi.
\end{equation*}
This gives the stiffness matrix
\begin{equation}
  \label{eq:ElementStiffness_timo}
  \bK_e
  =
  \frac{E_e I_e}{(1+\Phi_e) l_e^3}
  \mat[4]{
    12 	& 6l_e 		   & -12   & 6l_e \\
    6l_e & l_e^2 (4+\Phi_e) & -6l_e & l_e^2 (2-\phi_e) \\
    -12	 & -6l_e	   & 12	   & -6l_e \\
    6l_e & l_e^2 (2-\Phi_e) & -6l_e & l_e^2 (4+\Phi_e)
  }.
\end{equation}
The kinetic energy is
\begin{align*}
  T_e
  &= \frac12 \int_0^{l_e} \left(
    \rho_e A_e(\xi) \dot u_e(\xi, t)^2 +
    \rho_e I_e \dot\psi_e(\xi, t)^2 \right) d\xi \\
  &= \frac12 \int_0^{l_e} \left(
    \rho_e A_e(\xi) \dot u_e(\xi, t)^2 +
    \rho_e I_e [\dot\beta_e(\xi, t) + \dot u_e'(\xi, t)]^2 \right) d\xi.
\end{align*}
Again we write the approximate kinetic energy as
\begin{equation*}
  T_e(t) = \frac12 \dot q_e(t)^T \bM_e \dot q_e(t),
\end{equation*}
where the elements of the mass matrix are now
\begin{align*}
  m_{e,ij}
  &= \rho_e A_e \int_0^{l_e} N_{ei}(\xi) N_{ej}(\xi) \, d\xi + {} \\
  &\quad\; \rho_e I_e \int_0^{l_e} \left(
    \frac{\Phi_e l_e^2}{12} N_{ei}'''(\xi) + N_{ei}'(\xi) \right) \left(
    \frac{\Phi_e l_e^2}{12} N_{ej}'''(\xi) + N_{ej}'(\xi) \right) d\xi.
\end{align*}
Both integrals can be evaluated analytically, yielding
\begin{equation}
\label{eq:mass_matrix_timo}
\begin{aligned}
  \bM_e
  &= \frac{\rho_e A_e l_e}{840 (1+\Phi_e)^2}
    \mat[4]{
    m_1 & m_2 & m_3 & m_4 \\
    m_2 & m_5 & -m_4 & m_6 \\
    m_3 & -m_4 & m_1 & -m_2 \\
    m_4 & m_6 & -m_2 & m_5
    } + {} \\
  &\quad\; \frac{\rho_e I_e}{30 (1+\Phi_e)^2 l_e}
    \mat[4]{
    m_7 & m_8 & -m_7 & m_8 \\
    m_8 & m_9 & -m_8 & m_{10} \\
    -m_7 & -m_8 & m_7 & -m_8 \\
    m_8 & m_{10} & -m_8 & m_9
    },
\end{aligned}
\end{equation}
where
\begin{align*}
  m_1 &= 312 + 588 \Phi_e + 280 \Phi_e^2,
  & m_6 &= -(6 + 14 \Phi_e + 7 \Phi_e^2) l_e^2, \\
  m_2 &= (44 + 77 \Phi_e + 35 \Phi_e^2) l_e, & m_7 &= 36, \\
  m_3 &= 108 + 252 \Phi_e + 140 \Phi_e^2, & m_8 &= (3 - 15 \Phi_e) l_e, \\
  m_4 &= -(26 + 63 \Phi_e + 35 \Phi_e^2) l_e,
  & m_9 &= (4 + 5 \Phi_e + 10 \Phi_e^2) l_e^2, \\
  m_5 &= (8 + 14 \Phi_e + 7 \Phi_e^2) l_e^2,
  & m_{10} &= (-1 - 5 \Phi_e + 5 \Phi_e^2) l_e^2.
\end{align*}
The second matrix represents the effect of rotational inertia.
If shear effects and rotational inertia are ignored,
\eqref{eq:mass_matrix_timo} reduces to the Euler-Bernoulli mass matrix.

The mass matrix and stiffness matrix of the entire beam
are assembled in the same way as before,
and all mass matrices and stiffenss matrices are symmetric positive definite.

Note that the Young's modulus $E_e$
now appears not only in the stiffness matrix
but via $\Phi_e$ from \eqref{eq:shear-angle} also in the mass matrix.
Thus all these quantities may depend on the damage parameter vector $\bx$:
$E_e(\bx)$, $\Phi_e(\bx)$, $\bK_e(\bx)$, $\bM_e(x)$,
and finally $\bK(\bx)$ and $\bM(x)$.

\subsection{Damage Model}
\label{sec:model:damage}

\begin{figure}
  \centering
  \ximg[width=50mm,trim=0 8 0 24,clip]{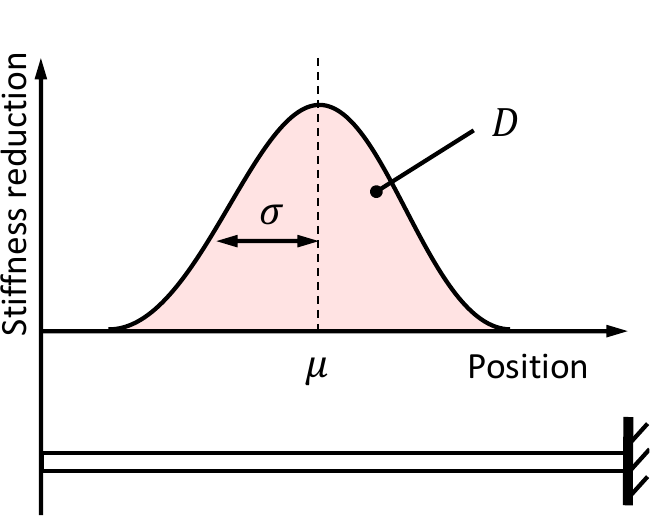}
  \caption{Damage distribution function.}
  \label{fig:damage-distribution}
\end{figure}

To model the damage we choose a univariate Gaussian distribution
described by three parameters $\bx \in \R^3$,
\begin{align*}
  \ddf{s}{\bx}
  &= \frac{D}{\sigma \sqrt{2 \pi}}
    \exp \left( -\frac{(s - \mu)^2}{2 \sigma^2} \right),
  &
    \bx &= (D, \mu, \sigma).
\end{align*}
Here $D$ is the damage severity (the total weight of the distribution),
$\mu$ is the center of the damage location,
and $\sigma$ is the standard deviation of the damage distribution,
a measure of the damage extent,
see \cref{fig:damage-distribution}.
In earlier work \cite{Wolniak23},
we have also experimented with
uniform distributions as damage models,
and the Gaussian distributions appear to be
the better match for most scenarios.
The parameters of primary interest
in typical structural health monitoring situations
are the location $\mu$ and the severity $D$
since the goal is to detect and locate small damage scenarios.
The extent $\sigma$ is basically a technical parameter which, however,
provides information on the accuracy of the location.
The severity has been chosen as the total weight
rather than the maximal value of $\ddf{}{}$,
\ie, $D / (\sigma \sqrt{2 \pi})$,
because this reduces the sensitivity of the model response
(eigenvalues and eigenfunctions)
with respect to the less relevant parameter $\sigma$.

The damage distribution function is used
to scale the initial stiffness of each beam element $e$.
This is accomplished by a scaling factor $\theta_e(\bx)$
that enters into \eqref{eq:ElementStiffness_Damaged}
and similarly into \eqref{eq:ElementStiffness_timo}
and via \eqref{eq:shear-angle} into \eqref{eq:mass_matrix_timo},
\begin{equation}
  E_e(\bx) = \theta_e(\bx) E_e^0,
\end{equation}
where $E_e^0$ is Young's modulus of beam element $e$ in the undamaged state.

Note that it is also possible to set up a more detailed model
with stiffness values for axial, shear, bending and torsional strains,
each depending individually on the damage.
However, such a model would yield an excessive number of damage parameters
that could not be determined efficiently in practice.

The scaling factors $\theta_e(\bx)$ are calculated using the
cumulative distribution function $\cpdf{}{}$ associated with $\ddf{}{}$,
\begin{equation*}
  \cpdf{s}{\bx} = \int_{-\oo}^s \ddf{r}{\bx} \, dr.
\end{equation*}
where $s$ is taken as a node position along the structure.
Following this approach, the stiffness scaling factors $\theta_e$
are defined as
\begin{equation}
  \label{eq:theta}
  \theta_e(\bx)
  =
  1 - L \frac{\cpdf{s_e}{\bx} - \cpdf{s_{e-1}}{\bx}}{l_e}.
\end{equation}
This results in a uniform distribution of the damage on each element.
In other words, it converts the smooth
distribution $\ddf{}{}$  (truncated to $[0, L]$)
into an element-wise constant distribution with the same weight,
$\int_0^L \ddf{r}{\bx} \, dr \lessapprox D$.

\subsection{Optimization Model}
\label{sec:model:opt}

The damage location problem has the basic structure of a
nonlinear inequality-constrained bi-objective optimization problem,
\begin{equation}
  \label{eq:opt}
  \min_{\bx \in \R^3} \ f(\bx)
  \qstq \bx \in [\bx^-, \bx^+],
  \quad c_e(\bx) \ge 0, \ e \in \setE.
\end{equation}
Here both components of $f$ are least-squares error measures,
see \cref{sec:model:lsq},
and in addition to component-wise lower and upper bounds $\bx^-$ and $\bx^+$
we have a nonlinear constraint at each element $e \in \setE$.
The functions $c_e\: \R^3 \to \R$ are smooth
whereas $f\: \R^3 \to \R^2$ is not differentiable
as it involves eigenvalues of the mechanical structure
which merely depend continuously on the parameter vector $\bx$.
We therefore combine the nonlinear constraints
into a single equivalent constraint,
\begin{equation}
  \label{eq:c}
  c(\bx) \define \min_{e \in \setE} c_e(\bx) \ge 0,
\end{equation}
which introduces artificial nonsmoothness but does not destroy
any structure that our pattern search algorithm could otherwise exploit.
In fact, although pattern search algorithms
cannot handle nonlinear equality constraints well,
nonlinear inequalities are properly incorporated
via infinite penalties at infeasible points \cite{coello2002}.
Thus, we finally replace \eqref{eq:opt} with
\begin{equation}
  \label{eq:opt-penalty}
  \min_{\bx \in \R^3} \ \_f(\bx)
  \qstq \bx \in [\bx^-, \bx^+],
  \qquad
  \_f(\bx)
  =
  \begin{cases}
    f(\bx), & c(\bx) \ge 0, \\
    \tbinom{+\oo}{+\oo}, & c(\bx) < 0.
  \end{cases}
\end{equation}

The specific bounds in our case are chosen as follows:
\begin{itemize}
\item $\mu \in [0, L]$: the damage center is on the structure;
\item $\sigma \in [0, L]$: the damage extent is in a reasonable range;
\item $D \in [0, D_{\max}]$ with $D_{\max} \le 1$:
  the damage severity is in a reasonable range.
\end{itemize}
The upper bound $L$ on $\sigma$ is artificial.
Since lower and upper bounds are required to define the search area
of our pattern search algorithm,
we choose a relatively large value to reduce the risk of missing any optima.
The choice of $D_{\max}$ is problem-specific:
if the damage is already known to be smaller,
a value $D_{\max} < 1$ reduces the range
that is to be searched by the algorithm.
The value $D_{\max} = 1$ corresponds to a ``full'' damage
that extends over all beam elements:
then $\sigma$ is at its upper bound
and all stiffness factors $\theta_e$ are close to zero.
This extreme case is very unlikely in practice
where both $D$ and $\sigma$ are typically small.

The specific nonlinear constraints enforce positive values
of the stiffness scaling factors,
\begin{equation}
  \label{eq:theta_constraint}
  c_e(\bx) \define \theta_e(\bx) - \theta\tsb{min}
  \qtextq{where} \theta\tsb{min} > 0.
\end{equation}
Positivity is not guaranteed by \eqref{eq:theta},
and in fact $\theta_e(\bx)$ can become negative for small values of $\sigma$,
which in turn may lead to an indefinite stiffness matrix $\bK(\bx)$.
The exact choice of $\theta\tsb{min}$ is application-specific.

In realistic cases the box constraints remain usually inactive,
with the exceptions of $D = 0$ (no damage)
and $\mu = 0$ or $\mu = L$ (damage at an end of the structure).
As to the nonlinear constraint \eqref{eq:c},
we observe no significant influence on the optimal solutions
for the relatively mild damage cases regarded in this work.
It rather becomes relevant only at the initial iterations
of our optimization algorithm, when the sampling pattern
is not yet focused close to a global optimum.
However, the constraint appears to help guiding the algorithm
away from unreasonable areas in the damage parameter space and,
thus, aids in increasing the convergence speed.

\subsubsection{Eigenvalue Problem}
\label{sec:model:eigenvalue}

We assume that an external excitation $\bm{p}(t)$ of the mechanical structure
generates only a small deviation $\bu(t)$ from some equilibrium state,
so that a standard linearly elastic damped oscillation
appropriately describes the dynamics,
\begin{equation}
  \label{eq:dynamic_mechanical_system}
  \bM(\bx) \ddot\bu(t) + \bB \dot\bu(t) + \bK(\bx) \bu(t) = \bm{p}(t).
\end{equation}
The modal parameters of the mechanical model are then obtained
from the generalized parameter-dependent eigenvalue problem
\begin{equation}
  \label{eq:eigenvalue}
  [\bK(\bx) - \lambda_k(\bx)^2 \bM(\bx)] \bu_k(\bx) = \bm0.
\end{equation}
Here $\lambda_k(\bx)$ are the eigenvalues,
and $\bu_k(\bx) \in \R^N$ are the corresponding eigenvectors,
\ie, spatially discretized eigenfunctions,
which are called \emph{mode shapes} in the engineering context.
The chosen formulation neglects damping and requires
that there is no external driving force.
Depending on the specific mechanical system, multiple eigenvalues may occur,
for instance when considering
bending modes of rotationally symmetric structures.
The oscillation frequencies corresponding to the eigenvalues are
\begin{equation}
  f_k(\bx) = \frac{\lambda_k(\bx)}{2 \pi}.
\end{equation}
It follows from \eqref{eq:eigenvalue} that
the eigenvalues $\lambda_k$ depend continuously
on the matrices $\bK$ and $\bM$ but in general not differentiably.
Thus, even though $\bK(\bx)$ and $\bM(\bx)$ are smooth,
$\lambda_k(\bx)$ is generally nonsmooth.

\subsubsection{Least-Squares Error Measures.}
\label{sec:model:lsq}

The measured data are time series of acceleration measurements
obtained from several sensors that are placed along the structure.
The sensor positions correspond to a subset of $m$ components
of the mode shapes $\bu_k(\bx)$,
which are selected by multiplication with a
gather matrix $\bS \in \set{0, 1}^{m \x N}$.
In addition, a normalization using the $\ell_2$-norm is carried out to obtain
\begin{equation}
  \label{eq:discrete-shape}
  \bPhi_k(\bx) = \dfrac{\bS \bu_k(\bx)}{\norm{\bS \bu_k(\bx)}} \in \R^m.
\end{equation}
The normalized mode shapes $\bPhi_k(\bx)$ can be directly compared
to the mode shapes identified from measurement data.
The identification of mode shapes $\bPhi_k$ and eigenfrequencies $f_k$
from acceleration measurements is achieved
by the Bayesian Operational Modal Analysis (BayOMA) \cite{Au2012}.

In the model updating framework,
four separate states and associated data sets
of the monitored structure are involved.
These states are generated by combining measured and simulated systems
with reference and analysis states.
The reference and analysis states correspond to
healthy and damaged structures, respectively.
The four combinations are given in \cref{tab:updating_states},
where the states are labeled with
M (measured), S (simulated), 0 (reference) and 1 (analysis).
For damage location, the states M0, M1 as well as S0 are constant
and only S1 is variable and subject to model updating.
\begin{table}
  \centering
  \caption{States of the structure considered in model updating.}
  \begin{tabular}{ccc}
    \toprule
    & Measurement & Simulation \\
    \midrule
    Reference (healthy) & M0 & S0 \\
    Analysis (damaged)  & M1 & S1 \\
    \bottomrule
  \end{tabular}
  \label{tab:updating_states}
\end{table}
The following two types of errors need to be considered
when setting up the least-squares objective:
\begin{enumerate}
\item Modeling errors in the simulation model (epistemic uncertainty)
  lead to different mechanical properties of measurement and simulation
  even in the reference state.
  Although these errors could in principle be reduced
  by more accurate modeling,
  a significant mismatch between measurement and simulation
  remains in practice.
  We refer to this modeling error in the reference state as
  ``error of order 0''.
\item The sensitivities of the real mechanical system and its simulation
  with respect to damage are roughly equal but not exactly.
  We refer to this modeling error in the sensitivity with respect to damage as
  ``error of order 1''.
\end{enumerate}
To eliminate the error of order 0
from the mismatch between simulation and measurement,
we define the mismatch in terms of the (normalized) differences
between the respective damaged and healthy states.
Then, using eigenfrequencies as well as mode shapes
to utilize the full information contained in these modal parameters,
we obtain the objective function
$f(\bx) = (\veps_f(\bx), \veps_m(\bx))$
in terms of two least-squares errors:
\begin{equation}
  \label{eqn:obj_fun}
  \begin{aligned}
    \veps_f(\bx)^2
    &=
    \sum_{k \in \setM} \left(
      \dfrac{f\sbS[k]1(\bx) - f\sbS[k]0}{f\sbS[k]0} -
      \dfrac{f\sbM[k]1 - f\sbM[k]0}{f\sbM[k]0}
    \right)^2,
    \\
    \veps_m(\bx)^2
    &=
    \sum_{k \in \setM} \norm{
      (\bPhi\sbS[k]1(\bx) - \bPhi\sbS[k]0) -
      (\bPhi\sbM[k]1      - \bPhi\sbM[k]0)
    }^2.
  \end{aligned}
\end{equation}
Here $\setM$ is the set of modes taken into consideration,
$\veps_f$ is the eigenfrequency error
and $\veps_m$ is the mode shape error.
As mentioned, the damage parameters only influence
the simulation results of the damaged case S1
while all other terms of \eqref{eqn:obj_fun} are constant.
All mode shapes in $\veps_m$ are normalized
according to \eqref{eq:discrete-shape},
similar to the approach of the
enhanced COMAC metric \cite{hunt1992application}
which has been proposed for structural monitoring.

%%% Local Variables:
%%% mode: latex
%%% TeX-master: t
%%% End:

\section{Optimization Algorithm}
\label{sec:algorithm}
In this section, we present the Multi-Objective Global Pattern Search Algorithm
for approximating solutions to multi-objective optimization problems
such as the bi-objective damage location problem given in \cref{sec:model:opt}.

\subsection{Multi-objective optimization}

Given a vector-valued function
\begin{equation*}
  f = (f_1, \dots, f_m)\: \Omega \to \R^m \cup \set{+\oo}^m
\end{equation*}
with $m \ge 2$ component functions
$f_1, \dots, f_m\: \Omega \to \R \cup \set{+\oo}$
and a nonempty feasible set $\Omega \subseteq \R^n$,
our general multi-objective optimization problem is given as
\begin{equation}
  \label{multiobjective_problem}
  \tag{MOP}
  \begin{cases}
    f(\bx) = (f_1(\bx), \dots, f_m(\bx)) \to \min, \\
    \bx\in \Omega.
  \end{cases}
\end{equation}
A feasible point $\bx \in \Omega$ is a Pareto-efficient (Pareto-optimal)
solution of the problem \eqref{multiobjective_problem}
if there is no $\_\bx \in \Omega$ such that
\begin{align*}
  f_i(\_\bx) & \le f_i(\bx) \qtext{for all } i = 1, \dots, m, \\
  f_j(\_\bx) & < f_j(\bx) \qtext{for some } j \in \set{1, \dots, m}.
\end{align*}
In other words, $\bx \in \Omega$ is a Pareto-efficient solution
if $f(\bx) \notin \set{+\oo}^m$ and
\begin{equation*}
  f[\Omega] \cap (\set{f(\bx)} - \R^m_+) = \set{f(\bx)},
\end{equation*}
where $f[\Omega] \define \defset{f(\_\bx)}{\_\bx \in \Omega}$
is the image of $\Omega$ under $f$
and $\R^m_+$ is the nonnegative orthant
(the so-called natural ordering cone) in $\R^m$.
The set of Pareto-efficient solutions is denoted by $\Eff(\Omega \mid f)$
and the set of non-dominated image points
(the so-called Pareto front or Pareto frontier)
is given by $\ND(f[\Omega])$, where
\begin{equation*}
  \ND(A) \define \defset{a \in A}{a \notin A \without{a} + \R^m_+}
\end{equation*}
for any set $A \subseteq \R^m \cup \set{+\oo}^m$.
Obviously, $f[\Eff(\Omega \mid f)] = \ND(f[\Omega])$.
The efficiency concept in multi-objective optimization dates back
to the works by Edgeworth \cite{Edgeworth1881} and Pareto \cite{Pareto1896}.
Some overview on multi-objective optimization can be found
in the books by Ehrgott \cite{Ehrgott2005} and Jahn \cite{Jahn11}
or in the survey by Eichfelder \cite{Eichfelder2021a}.

\begin{figure}
  \centering
  \ximg[height=36mm]{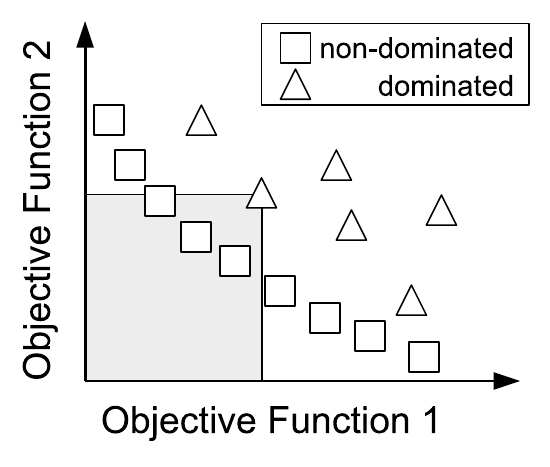}
  \caption{Sketch of objective value space
    in multi-objective optimization and Pareto dominance.}
  \label{fig:non_dom}
\end{figure}

In \cref{fig:non_dom}, the points marked by squares
constitute the set of non-dominated image points,
while the points marked by triangles are dominated by several other points
and are therefore not on the Pareto front.

\subsection{Penalty approach for the nonlinear constraints}

In many practical problems,
such as our model updating problem given in \cref{sec:model:opt},
$\Omega$ is given by a set of the type
\begin{equation}
  \label{eq:constraint_set_g}
  \defset{\bx \in [\bx^-, \bx^+]}{c_1(\bx) \ge 0, \dots, c_l(\bx) \ge 0},
\end{equation}
where $c_1, \dots, c_l\: [\bx^-, \bx^+] \to \R$, $l \in \N$,
are constraint functions,
and $\bx^-$ and $\bx^+$ are component-wise lower and upper bounds.
Penalty methods in multi-objective optimization
aim to transform the problem \eqref{multiobjective_problem}
with feasible set given in \eqref{eq:constraint_set_g}
into a corresponding problem of type \eqref{multiobjective_problem}
involving a single box constraint $\bx \in [\bx^-, \bx^+]$.
Consider a penalty function
$\varphi\: [\bx^-, \bx^+] \to \R_+ \cup \set{+\oo}$
with the property that for any $\bx \in [\bx^-, \bx^+]$,
\begin{equation*}
  \varphi(\bx) = 0 \quad \iff \quad c_1(\bx) \ge 0, \dots, c_l(\bx) \ge 0.
\end{equation*}
For instance, in numerical penalty methods
the function $\varphi$ is often defined by
\begin{equation*}
  \varphi(\bx) \define -\min\set{0, c_1(\bx), \dots, c_l(\bx)}
  \qtextq{or by}
  \varphi(\bx) \define \sum_{i = 1}^l \min\set{0, c_i(\bx)}^2.
\end{equation*}
For given (sufficiently large) penalty parameters
$\alpha_1, \dots, \alpha_m > 0$
we consider the penalized multi-objective optimization problem
\begin{equation}
  \label{multiobjective_problem_penalised}
  \tag{MOP$\tsb{pen}$}
  \begin{cases}
    (f + \alpha \varphi)(\bx)
    =
    (f_1(\bx) + \alpha_1 \varphi(\bx),
    \dots,
    f_m(\bx) + \alpha_m \varphi(\bx)) \to \min, \\
    \bx\in [\bx^-, \bx^+].
  \end{cases}
\end{equation}
In the literature, various relationships are established between the problems
\eqref{multiobjective_problem} and \eqref{multiobjective_problem_penalised};
see, \eg, Ye \cite{Ye2012} and G\"unther et al.\ \cite{GuentherEtAl2023}.
From a theoretical point of view, one can also apply
the method of infinite penalty (extreme barrier approach), which is based on
\begin{equation*}
  \varphi(\bx)
  \define
  \begin{cases}
    0 & \text{for } \min\set{c_1(\bx), \dots, c_l(\bx)} \ge 0, \\
    +\oo & \text{otherwise},
  \end{cases}
\end{equation*}
and $\alpha_1 = \dots = \alpha_m  = 1$.

Based on the above considerations, in the following sections
we consider the problem \eqref{multiobjective_problem}
with a feasible set $\Omega$ given by a single box constraint,
$\Omega \define [\bx^-, \bx^+]$,
and we employ the method of infinite penalty.

\subsection{Basic ideas of the pattern search approach}
\label{sec:basic_ideas_algorithm}

\begin{figure}
  \centering
  \ximg[height=45mm]{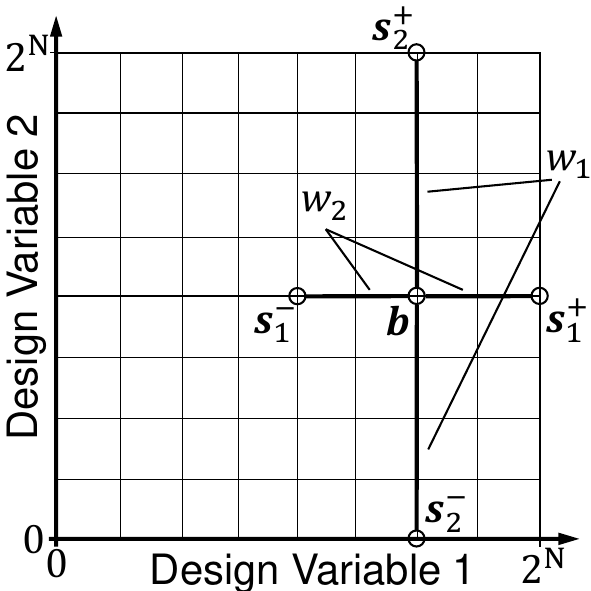}
  \caption{Sampling pattern of our Multi-Objective
    Global Pattern Search Algorithm in two dimensions.}
  \label{fig:sampling}
\end{figure}

In our pattern search approach, we create a pattern as usual
(for such a type of algorithm)
by performing a one-at-a-time permutation
around a base vector $\bb = (\bb_1, \dots, \bb_n)$
using a step width vector $\bw = (\bw_1, \dots, \bw_n)$.
New pattern sampling points are then defined
for each $i = 1, \dots, n$ by the vectors
\begin{align*}
  s_i^-(\bb) \define \bb - \bw_i e^i \qtextq{and}
  s_i^+(\bb) \define \bb + \bw_i e^i,
\end{align*}
where $e^i$ is the $i$-th canonical basis vector in $\R^n$.
To represent the resulting grid in floating-point arithmetic,
we need to avoid unnecessary rounding errors.
For instance, the equality $\bb^1 + \bw_i e^i = \bb^2 - \bw_i e^i$
with base vectors $\bb^1$ and $\bb^2$
might hold in exact arithmetic but not in floating-point arithmetic.
For this reason, we introduce the box $[0, 2^N]^n$
and basis vectors with integer coordinates
to obtain the following transformed multi-objective optimization problem:
\begin{equation}
  \label{multiobjective_problem_g}
  \tag{MOP$_g$}
  \begin{cases}
    g(s) = (f_1(h(s)), \dots, f_m(h(s))) \to \min\\
    s \in [0, 2^N]^n.
  \end{cases}
\end{equation}
Here $h\: [0, 2^N]^n \to [\bx^-, \bx^+]$ is a bijective linear function
defined as
\begin{equation*}
  h_i(s) \define \bx^-_i + 2^{-N} s_i (\bx^+_i - \bx_i^-)
  \qtext{for all  } s = (s_1, \dots, s_n) \in [0, 2^N]^n, \, i = 1, \dots, n,
\end{equation*}
and $g\: [0, 2^N]^n \to \R^m \cup \set{+\oo}^m$ is defined as
\begin{equation*}
  g(s) \define f(h(s)) \qtext{for all  } s \in [0, 2^N]^n.
\end{equation*}
It is easy to check that $f[[\bx^-, \bx^+]] = g[[0, 2^N]^n]$, and so
\begin{equation*}
  \ND(f[[\bx^-, \bx^+]]) = \ND(g[[0, 2^N]^n]).
\end{equation*}
Moreover, using the bijection property of $h$ we get
\begin{equation*}
  \Eff([\bx^-, \bx^+] \mid f) = h[\Eff([0, 2^N]^n \mid g)].
\end{equation*}
When applying a pattern search scheme
to the problem \eqref{multiobjective_problem_g}, the set
\begin{equation*}
  G_\Z \define ([0, 2^N] \cap \Z)^n,
\end{equation*}
which consists of all integer grid points in the box $[0, 2^N]^n$,
is of particular interest.
In our approach we will use step widths
$\bw_1, \dots, \bw_n \in \set{2^0, 2^1, \dots, 2^{N-1}}$,
and so if $\bb \in G_\Z$ and $s_i^-(\bb), s_i^+(\bb) \in [0, 2^N]^n$,
then $s_i^-(\bb), s_i^+(\bb) \in G_\Z$.
Roughly speaking, our Multi-Objective Global Pattern Search Algorithm
is an iterative procedure based on the following tasks:
\begin{itemize}
\item Define a set of base vectors
  (to start, take the point $\bb = 2^{N-1} (1, \dots, 1)$,
  which is the center of the box $[0, 2^N]^n$,
  and the step width vector $\bw = \bb$).
\item Sample new points in $G_\Z$
  (with the aid of $s_i^-(\bb)$ and $s_i^+(\bb)$)
  around all base vectors.
\item Apply a non-dominated sorting procedure
  to all current base vectors and new sampling points
  to obtain a new set of base vectors for the next iteration.
\item Halve the maximum step width of one variable
  if the set of base vectors does not change in two iterations.
\end{itemize}
The algorithm stops if the maximum of the step widths is equal to one.
By taking a look at the above approach (in particular, the definition of $h$),
one can see that a larger parameter $N$ leads to a larger integer grid $G_\Z$
and so to a finer grid in the original box $h[[0, 2^N]^n] = [\bx^-, \bx^+]$.
The integer grid $G_\Z$,
the pattern generation scheme in the box $[0, 2^N]^n$,
and the resulting sampling points
$s_1^-(\bb), s_1^+(\bb), s_2^-(\bb), s_2^+(\bb)$
around a base vector $\bb$
are illustrated in \cref{fig:sampling} for a two-dimensional example.

The process of identifying the non-dominated set $\ND(A)$
of a finite set of points, such as
$A = \set{a^1, \dots, a^p} \subseteq f[\Omega] \cup \set{+\oo}^m$,
is called non-dominated sorting in the literature \cite{kalyanmoy2001multi}.
In our Multi-Objective Global Pattern Search Algorithm,
we use the Graef-Younes method with presorting
(see G\"unther and Popovici \cite{Guenther18},
and also Jahn \cite[Sec. 12.4]{Jahn11})
for the non-dominated sorting phase,
which we recall in the following part of the section.

Consider a vector $\lambda = (\lambda_1, \dots, \lambda_m) \in \R^m$,
$\lambda_i > 0$ for all $i = 1, \dots, m$,
and a function $\varphi_\lambda\: A \to \R$ defined by
\begin{equation*}
  \varphi_\lambda(a) \define \lambda^T a \qtext{for all } a \in A.
\end{equation*}
In multi-objective optimization,
the strongly $\R^m_+$-increasing function $\varphi_\lambda$
is used in the weighted sum (linear) scalarization method
(see Jahn \cite{Jahn11}).
For the given choice of $\lambda$ one has
\begin{equation*}
  \arg\min_{a \in A} \varphi_\lambda(a) \subseteq \ND(A).
\end{equation*}

The Graef-Younes method \cite[Sec. 12.4]{Jahn11}
given in the following \cref{alg:GYmethod}
produces a set $B \subseteq A$ with the property $\ND(A) \subseteq B$.
\begin{algorithm}
  \caption{Graef-Younes reduction method (function \GYRM)}
  \label{alg:GYmethod}
  \algsetup{indent=2em}
  \begin{algorithmic}
    \STATE Input: $A = \set{a^1, \dots, a^p}$
    \STATE $B \gets \set{a^1}$
    \FOR{$k \gets 2$ \TO $p$}
    \IF{$a^k \notin B + \R^m_+$}
    \STATE $B \gets B \cup \set{a^k}$
    \ENDIF
    \ENDFOR
    \STATE Output: $B \subseteq A$ with $\ND(A) \subseteq B$
  \end{algorithmic}
\end{algorithm}

In our proposed algorithm, we use the Graef-Younes method with presorting
to generate the set $B = \ND(A)$, as shown in \cref{alg:preSortGYmethod}
(see G\"unther and Popovici \cite{Guenther18}).
\begin{algorithm}
  \caption{Graef-Younes method with presorting procedure}
  \label{alg:preSortGYmethod}
  \algsetup{indent=2em}
  \begin{algorithmic}
    \STATE Input: $A = \set{a^1, \dots, a^p}$
    \STATE Compute an enumeration $A = \set{a^{j_1}, \dots, a^{j_p}}$ such that
    \begin{equation*}
      \varphi_\lambda(a^{j_1}) \le
      \varphi_\lambda(a^{j_2}) \le \dots \le \varphi_\lambda(a^{j_p}).
    \end{equation*}
    \vspace*{-2ex}
    \STATE $B \gets \GYRM(\set{a^{j_1}, a^{j_2}, \dots, a^{j_p}})$
    \STATE Output: $B = \ND(A)$
  \end{algorithmic}
\end{algorithm}

\Cref{alg:preSortGYmethod} is an output-sensitive method, \ie,
the runtime depends (in addition to the size of the input $A$)
on the size of the output $B$.
More precisely, the worst-case complexity according to \cite{Guenther18}
can be specified as follows:
\begin{equation*}
  \bigO(|A| \log|A| + |A| \, |B|).
\end{equation*}
If $|B|$ is very small (\eg, $|B| \le \log|A|$),
then \cref{alg:preSortGYmethod} is very efficient.
The calculation of the values
$\varphi_\lambda(a^{j_i})$, $i \in \set{1, \dots, p}$,
can be efficiently performed by parallelization.

Besides the Graef-Younes method with presorting,
also some other non-dominated sorting methods are known in the literature;
see, \eg, \cite{Kung1975} and \cite{deb2000}.

\begin{figure}
  \centering
  \ximg[height=38mm]{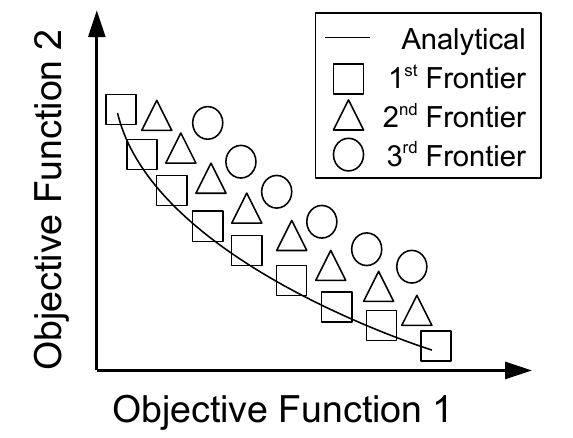}
  \caption{Illustration of higher level Pareto fronts
    resulting from repeated non-dominated sorting.}
  \label{fig:2ndfront}
\end{figure}

For pattern search type algorithms it is known
that the concept of higher level Pareto fronts is useful
to achieve a better approximation of the set of efficient solutions
\cite{kalyanmoy2001multi}.
Note that the first level Pareto front is given by $B_1 \define \ND(A)$,
the second level Pareto front by $B_2 \define \ND(A \setminus B_1)$,
and generally the level $i$ Pareto front for each $i > 1$ by
$B_i \define \ND(A \setminus (B_1 \cup \dots \cup B_{i-1}))$.
\Cref{fig:2ndfront} illustrates this concept
and shows how a second level and third level front
emerge in parallel to the first level front.
These higher level fronts can be computed by \cref{alg:ParetoFrontiers}.

\begin{algorithm}
  \caption{Computing Pareto fronts $1$ to $k$}
  \label{alg:ParetoFrontiers}
  \algsetup{indent=2em}
  \begin{algorithmic}
    \STATE Input: $A = \set{a^1, \dots, a^p}$, $k \in \N$
    \STATE Compute an enumeration $A = \set{a^{j_1}, \dots, a^{j_p}}$ such that
    \begin{equation*}
      \varphi_\lambda(a^{j_1}) \le
      \varphi_\lambda(a^{j_2}) \le \dots \le \varphi_\lambda(a^{j_p}).
    \end{equation*}
    \vspace*{-2ex}
    \FOR{$i = 1$ \TO $k$}
    \STATE $B_i \gets \GYRM(A \setminus \bigcup_{j = 1}^{i-1} B_j)$
    \ENDFOR
    \STATE Output: $B_1, \dots, B_k$ represent Pareto fronts 1 to k
  \end{algorithmic}
\end{algorithm}
\Cref{alg:ParetoFrontiers} is also an output-sensitive method,
with a worst-case complexity of
\begin{equation*}
  \bigO\Bigl( |A| \log|A| + |A| \sum_{i = 1}^k |B_i| \Bigr).
\end{equation*}
The complexity of calculating all Pareto fronts is therefore given by
\begin{equation*}
  \bigO(|A| \log|A| + |A|^2).
\end{equation*}
Note that \cref{alg:ParetoFrontiers} is a novel alternative
to the method proposed in \cite{deb2000}
for computing all Pareto fronts with quadratic worst-case complexity.

\subsection{Formulation of the Multi-Objective Global Pattern Search Algorithm}

\begin{figure}
  \centering
  \ximg[height=40mm]{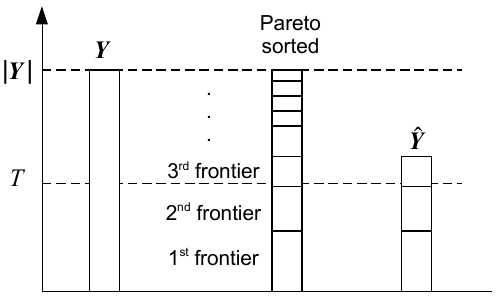}
  \caption{Illustration of the `Hall of Fame' updating scheme
    based on non-dominated sorting and the parameter $T$.}
  \label{fig:hof_update}
\end{figure}

In this section, we present our
Multi-Objective Global Pattern Search Algorithm in full detail.
First, as already described in the general scheme of the algorithm
in \cref{sec:basic_ideas_algorithm},
we use a non-dominated sorting procedure
at each iteration of the algorithm (loop \dots\ end loop).
In \cref{alg:points} (update of the `Hall of Fame')
we describe the details of the non-dominated sorting procedure
used to generate (first and higher level) Pareto fronts.
Starting from a set $\bY$ of sampling points,
the `Hall of Fame' procedure generates Pareto fronts
(starting with level one and then gradually increasing)
until at least $\min\set{T, |\bY|}$ elements
are in the union of the considered Pareto fronts.
This union $\^\bY$ of Pareto fronts is then used
in the Multi-Objective Global Pattern Search Algorithm
to define the next set of basis vectors of the pattern search scheme.
The update scheme for the `Hall of Fame' is shown schematically
in \cref{fig:hof_update}.

\begin{algorithm}
  \caption{Update the `Hall of Fame' (function \uHoF)}
  \label{alg:points}
  \algsetup{indent=2em}
  \begin{algorithmic}
    \STATE Input: The set $\bY$
    \STATE Compute an enumeration $\bY = \set{y^{j_1}, \dots, y^{j_{|\bY|}}}$
    such that
    \begin{equation*}
      \varphi_\lambda(y^{j_1}) \le
      \varphi_\lambda(y^{j_2}) \le \dots \le \varphi_\lambda(y^{j_{|\bY|}}).
    \end{equation*}
    \vspace*{-2ex}
    \STATE $\^\bY \gets \emptyset$
    \WHILE{$|\^\bY| < \min\set{T, |\bY|}$}
    \STATE $\^\bY \gets \^\bY \cup \GYRM(\bY \setminus \^\bY)$
    % \COMMENT{Add non-dominated set to hall of fame.}
    \ENDWHILE
    \STATE Output: The set $\^\bY$
  \end{algorithmic}
\end{algorithm}

Now we are ready to formulate the complete
Multi-Objective Global Pattern Search Algorithm (MOGPSA for short)
which is given in \cref{alg:mogps_new}.
\begin{algorithm}
  \caption{Multi-Objective Global Pattern Search}
  \label{alg:mogps_new}
  \algsetup{indent=2em}
  \begin{algorithmic}
    \STATE Input: Problem \eqref{multiobjective_problem}
    (respectively, \eqref{multiobjective_problem_g}), parameters $T, N \in \N$
    \STATE $\bw \gets 2^{N-1} (1, \dots, 1)$
    \STATE $\bB \gets \set{\bw}$
    \STATE $\bB\tsb{all} \gets \bB$
    \LOOP
    \STATE $S \gets \emptyset$
    \FOR{$\bb \in \bB$}
    \FOR{$i = 1$ \TO $n$}
    \IF{$s_i^+(\bb) \in [0, 2^N]^n \setminus \bB\tsb{all}$}
    \STATE $S \gets S \cup \set{s_i^+(\bb)}$
    \ENDIF
    \IF{$s_i^-(\bb) \in [0, 2^N]^n \setminus \bB\tsb{all}$}
    \STATE $S \gets S \cup \set{s_i^-(\bb)}$
    \ENDIF
    \ENDFOR
    \ENDFOR
    \STATE $\bB\tsb{all} \gets \bB\tsb{all} \cup S$
    \STATE $\bY \gets g[\bB \cup S]$
    \STATE $\^\bY \gets$ \uHoF($\bY$)
    \STATE $\^\bB \gets g^{-1}[\^\bY] \cap (\bB \cup S)$
    \IF{$\^\bB \ne \bB$}
    \STATE $\bB \gets \^\bB$
    \STATE \textbf{continue} loop
    \ENDIF
    \IF{$\max\set{\bw_1, \dots, \bw_n} = 1$}
    \STATE \textbf{break} loop
    \ENDIF
    \STATE Choose $j \in \set{1, \dots, n}$
    such that $\bw_j = \max\set{\bw_1, \dots, \bw_n}$
    \STATE $\bw_j \gets \bw_j / 2$
    \ENDLOOP
    \STATE$\^\bY \gets$ \GYRM($\^\bY$)
    \STATE $\^\bB \gets g^{-1}[\^\bY] \cap \^\bB$
    \STATE Output: The sets $h[\^\bB] \subseteq [\bx^-, \bx^+]$
    and $\^\bY \subseteq f[[\bx^-, \bx^+]]$
  \end{algorithmic}
\end{algorithm}
We would like to point out that \cref{alg:mogps_new} terminates.
Indeed, consider an arbitrary iteration $k$ of \cref{alg:mogps_new}
with underlying step widths $\bw_1, \dots, \bw_n$.
We claim that there is an iteration $k+l$ (with $l \in \N \cup \set{0}$)
so that $\^\bB = \bB$, hence (if $\max\set{\bw_1, \dots, \bw_n} > 1$)
the algorithm reaches the last loop line
in which the maximum step width of one variable is divided by two.
Of course, the condition $\^\bB = \bB$ can occur
in an iteration $k+l$ with $S \neq \emptyset$.
Otherwise, $\^\bB = \bB$ occurs in an iteration $k+l$ with $S = \emptyset$,
taking into account that there are only finitely many
integer grid points in $[0, 2^N]^n$.
Consequently, the step widths $\bw_1, \dots, \bw_n$
fulfill the condition $\max\set{\bw_1, \dots, \bw_n} = 1$
after a finite number of iterations.
The break command then ends the loop.

\subsection{Comments related to the algorithm}

First of all, we would like to reflect
about the choice of the algorithm parameter $T$.

Clearly, for $T \ge 1$ the set $\^\bY$
generated by the 'Hall of Fame' procedure
(as part of \cref{alg:mogps_new})
always contains the first level Pareto front $\ND(\bY)$.

If $T = 1$, then $\min\{T, |\bY|\} = 1$,
so the 'Hall of Fame' procedure generates nothing else
than the first level Pareto front $\ND(\bY)$.

If $T \ge |G_{\Z}|$ one has $\^\bB = \bB \cup S$
after the 'Hall of Fame' procedure.
Then, at the end of the algorithm (after the loop),
a non-dominated (first level) sorting is applied
to all points in $g[G_{\Z}]$, \ie, $\^\bY = \ND(g[G_{\Z}])$.
If $N$ has a sufficiently large value,
one obtains a good approximation of the box $h[[0,2^N]^n] = [\bx^-, \bx^+]$
and thus a good approximation of the set of efficient solutions
$\Eff([\bx^-, \bx^+] \mid f)$.

Based on the above observations (and some numerical experience),
one can expect that our algorithm achieves a better global behavior
with respect to the approximation of the set of efficient solutions
$\Eff([\bx^-, \bx^+] \mid f)$
if one increases $T$ (and chooses a sufficiently large value for $N$).
Of course, the drawback is that an increased parameter $T$
often leads to a higher numerical effort
for computing the sets $\^\bB$ and $\^\bY$.
In particular, this is the case
if the calculation of the objective function values
of the given sampling points is very costly,
as in our bi-objective damage location problem in \cref{sec:model:opt}.
It is clear that a parallel evaluation of the function values
of the given sampling points is useful,
especially for a larger $T$ with a higher number of sampling points.

In practice, when choosing $T$ (and $N$) for a specific problem,
the decision-maker has to gain experience with the \cref{alg:mogps_new}
with regard to the expected approximation quality and the time supply.

Considering that calculating the objective function values
of sampling points can be very costly in practice,
it is important in an effective implementation of \cref{alg:mogps_new}
to store the assignment between sampling points
and corresponding calculated image points (using a cache structure).
This means that one should avoid to recalculate
$g(\bb)$ for $\bb \in \bB \cup S$
(statement $\bY \leftarrow g[\bB \cup S]$ in \cref{alg:mogps_new}),
if $\bb$ is already known from previous iterations.
Using the assignment between sampling points and corresponding image points,
the statements
$\^\bB \leftarrow g^{-1}[\^\bY] \cap (\bB \cup S)$ and
$\^\bB \leftarrow g^{-1}[\^\bY] \cap \^\bB$,
which are given in \cref{alg:mogps_new},
can be effectively implemented as well.

In the following part of the section,
we relate our approach to other approaches known from the literature.

First of all, if we apply MOGPSA to the single-objective case
(\ie, $m = 1$ in the problem \eqref{multiobjective_problem}),
we basically obtain the (Single-Objective) Global Pattern Search Algorithm
(GPSA) proposed in \cite{gps}.
In this case, for each set $A \subseteq \R \cup \set{+\oo}$, we have
\begin{equation*}
  \ND(A)
  =
  \defset{a \in A}{a \notin A \without{a} + \R^m_+}
  =
  \defset{a \in A}{\nexists \bar a \in A\: \bar a < a}.
\end{equation*}
The non-dominated sorting procedure
according to the \cref{alg:preSortGYmethod}
is therefore basically a sorting of finitely many numbers
in $\R \cup \set{+\oo}$ and the selection of the smallest element.

MOGPSA belongs to the class of deterministic derivative-free
pattern search methods for approximating solutions
to multi-objective optimization problems,
especially problems with expensive objective functions.
Clearly, it shares the pattern search feature with other algorithms
known in (single and multi-objective) optimization; see, \eg,
\cite{Hooke1961,torczon1997,Dolan2003,abramson2005,custodio2011,
  custodio2012,gps} and the references therein.
According to \cite{custodio2012} our algorithm
belongs to the class of direct search methods.
Some well-known deterministic derivative-free methods
for multi-objective optimization problems are proposed in
\cite{custodio2011,custodio2012,evtushenko2014},
while non-deterministic (hybrid) methods are proposed in
\cite{deb2000,Jahn2006,custodio2011},
and many other publications related to evolutionary
(population-based) multi-objective optimization;
see \cite{custodio2012} for a survey on this topic.
Non-dominated sorting for computing first and higher level Pareto fronts
based on a new procedure with quadratic worst-case complexity
(see \cref{alg:ParetoFrontiers}) is a key feature of our algorithm,
similar to NSGA-II in \cite{deb2000} but in contrast to the
Direct Multisearch Method in \cite{custodio2011}
and MOSAST in \cite{Jahn2006,Limmer2012}
where only the first level Pareto front is used.
MOGPSA has only two input parameters $T$ and $N$
which have an influence on the quality of the approximation of solutions
of the multi-objective optimization problem \eqref{multiobjective_problem}.
Numerical experience with our \cref{alg:mogps_new} shows
that the set $\^\bY$ generated by the 'Hall of Fame' procedure
is often larger at late iterations,
since then many first-level (and possibly higher level)
non-dominated points are contained in $\^\bY$
(notice again that we always have $\ND(\bY) \subseteq \^\bY$
in the 'Hall of Fame' procedure),
which is a useful feature of our algorithm to obtain (in the best case)
a good approximation of the Pareto front
of the problem \eqref{multiobjective_problem}.
However, as noted above, the choice of $T$ is also important
in order to have sufficiently many points in $\^\bY$
generated by the 'Hall of Fame' procedure in early iterations.
In particular, this means that the number of (first or higher level)
non-dominated points considered at each iteration
of our \cref{alg:mogps_new} generally does not remain constant.
This is in contrast to certain evolutionary algorithms,
such as NSGA-II in \cite{deb2000},
where only a constant number of non-dominated points are used
at each iteration of such algorithms
(possibly not including complete Pareto fronts).
For instance, such a reduction in the number of (first or higher level)
non-dominated points is achieved in NSGA-II by using a crowding distance metric.

\section{Computational Results}
\label{sec:results}
We test the proposed optimization algorithm
and demonstrate its applicability
on two example problems: a cantilever beam and a girder mast.
In both cases we introduce a controllable artificial damage
by local modifications of the structural properties.
As our focus is on the practial application here,
we do not evaluate the algorithmic performance.

\subsection{Laboratory beam experiment}

\begin{figure}
  \centering
  \ximg[width=\linewidth,trim=8 160 6 42,clip]% Relative trim units.
  {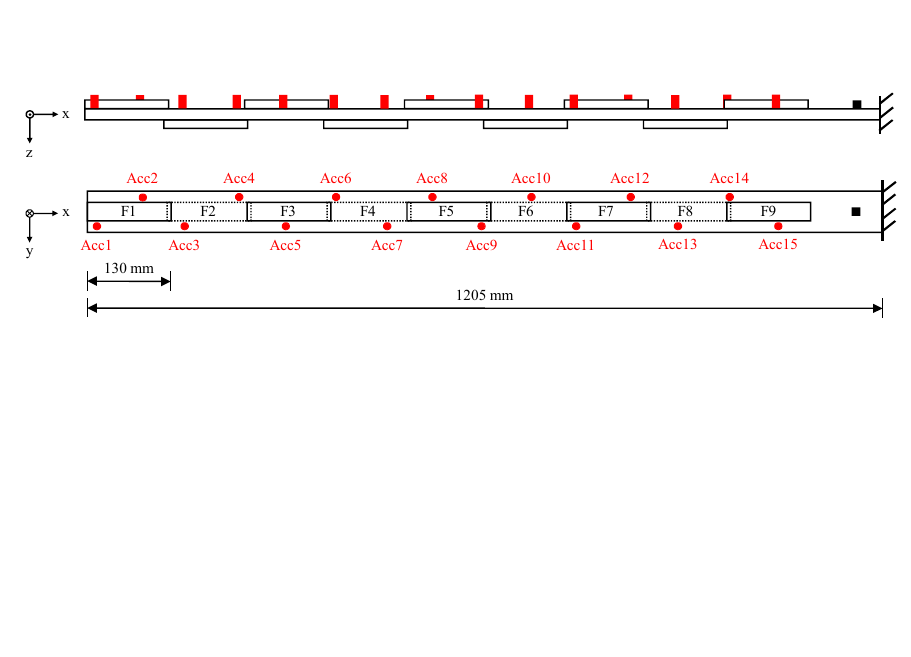}
  \caption{Schematic overview of the steel cantilever beam
    with nine screwed-on fishplates F$1$ to F$9$,
    15 sensor positions (red points),
    and the position of excitation (black square).}
  \label{fig:SketchBeam}
\end{figure}

\begin{figure}
  \centering
  \ximg[width=.98\linewidth]{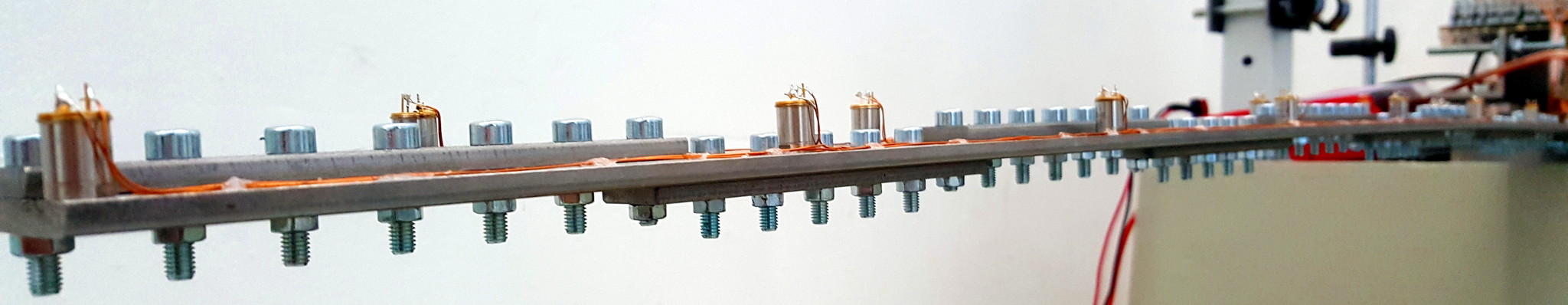}
  \caption{Photograph of the steel cantilever beam.}
  \label{fig:PicBeam}
\end{figure}

\begin{figure}
  \centering
  \ximg[angle=90, width=100mm]{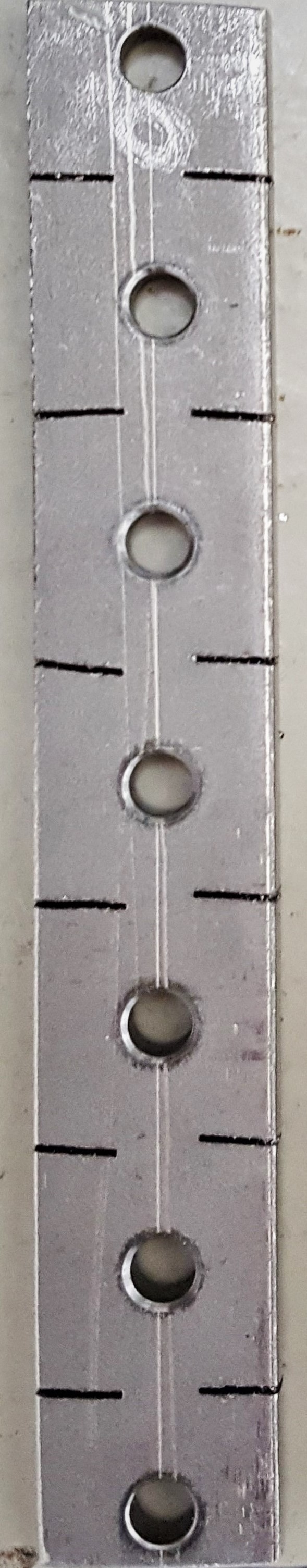}
  \caption{Photograph of reversible damage mechanism for cantilever beam,
    \ie, damaged fishplate.}
  \label{fig:damage_fishplate}
\end{figure}

As a first numerical example,
a laboratory cantilever beam structure is chosen.
Measurement data obtained from an experimental steel beam
with reversible damage mechanisms
presented in Wolniak et al.\ \cite{Wolniak23} is utilized.
\Cref{fig:SketchBeam,fig:PicBeam} show a schematic overview
as well as a photograph of the cantilever beam;
the artificial damage mechanism is shown in \cref{fig:damage_fishplate}.

\subsubsection{Experimental Setup}

Nine fishplates are mounted to the top and bottom of a central beam.
They can be replaced with damaged fishplates
to reduce the stiffness at various positions
(F1 to F9 in \cref{fig:SketchBeam}).
The fishplates are damaged by uniformly distributed sawcuts
on both sides of the plate, as illustrated in \cref{fig:damage_fishplate}.
Thus, replacing an undamaged fishplate with a damaged one
reduces the stiffness at known positions.
The damage leads to a decrease in the eigenfrequencies of the structure
as well as to alterations in the bending mode shapes.
In order to identify the modal parameters,
a contact-free electromagnetic shaker is positioned
close to the root of the beam (black square in \cref{fig:SketchBeam})
and used to excite the structure
with a broadband white noise excitation signal with up to \SI{250}{Hz}.
Fifteen IEPE acceleration sensors of the type MEAS 805-0005-01
are placed uniformly along the length of the beam
and record the resulting vibration
(Acc1 to Acc15 in \cref{fig:SketchBeam}).
The electrical signals provided by the sensors are recorded
with an HBK QuantumX MX1601B data acquisition system
using analog-to-digital conversion with 24-bit resolution
and a sampling rate of \SI{1200}{Hz}.
A detailed description of the experiment is given
in Wolniak et al. \cite{Wolniak23}.
The modal parameters $\bPhi_k$ and $f_k$ are identified from
time-domain measurements with a duration of one minute each
as described in \cref{sec:model:opt}.
For model updating, the identification results
obtained for the reference state
represent the measured healthy state M0,
as defined in \cref{tab:updating_states}.
The identification results of the nine damage scenarios
are considered as the analysis state M1.

The beam is uniformly discretized into 241 elements of length \SI{5}{mm}
which obey Ber\-noul\-li linear beam theory (\cref{sec:model:struct}).
As the cantilever beam has sections of varying thickness,
as well as additional sensor and screw masses,
the cross-sectional area $A_e$, the linear density $l_e \rho_e$
and the area moment of inertia $I_e$ of the beam elements
must be calculated accordingly.
The fishplates are mounted on either side of the beam
in an alternating pattern
with a small overlap between them along the length of the beam,
resulting in beam elements with a combined cross-sectional area
of the beam and one fishplate along with elements
with a combined cross-sectional area of the beam and two fishplates.
The varying cross sections are taken into account
when calculating the beam element's area moment of inertia.
This is done by adding either one or two times the fishplate's
moment of inertia on top of the beam's.
Finally the linear density is calculated,
accounting for varying cross sections.
In addition, the sensor masses are considered
by adding them in the form of an equivalent linear density
at the appropriate positions.
Screw masses are accounted for by a uniformly distributed
linear density on every element.

The finite element model is created and assembled in \Matlab
using the varying properties of the beam elements.
Mechanical parameters of this structure
are given in \cref{tab:cantilever_properties}.
Damage is modeled by reducing the element stiffness
according to the Gaussian damage distribution
with parameters $\bx = (D, \mu, \sigma)$
presented in \cref{sec:model:damage}.
\begin{table}
  \centering
  \caption{Properties of the steel cantilever beam.}
  \begin{tabular}{lr}
    \toprule
    Property & Value \\
    \midrule
    Elastic Modulus & \SI{127}{GPa} \\
    Beam Length  & \SI{1205}{mm} \\
    Beam Height  & \SI{5.15}{mm} \\
    Beam Width   & \SI{60}{mm} \\
    Fishplate Height  & \SI{4.85}{mm} \\
    Fishplate Width   & \SI{2}{mm} \\
    Density & \SI{7800}{kg/m^3}\\
    \bottomrule
    \end{tabular}
    \label{tab:cantilever_properties}
\end{table}

We consider nine damage scenarios obtained by replacing
exactly one of the fishplates F$i$ with a damaged one.
The Gaussian distribution model intentionally differs
from the uniform distribution used to damage the fishplates.
The damage distribution is parameterized according to the formulation
introduced in \cref{sec:model:struct} using the center of damage $\mu$,
the damage severity $D$ and the damage extent $\sigma$.
For each of the nine damage scenarios we know the center $\mu$
but we cannot compute the exact values of $\sigma$ and $D$.
The domain of optimization is given by the bounds
described in \cref{sec:model:damage},
with $D_{\max} = 0.3$ for the severity $D \in [0, D_{\max}]$
and with length $L = \SI{1205}{mm}$ for $\mu, \sigma \in [0, L]$.
The constraint for the stiffness scaling factors
is enforced below $\theta_\text{min}=0.15$.
These values are chosen such that the considered damage scenarios
are situated inside the domain.

\subsubsection{Computational Results}

\begin{figure}
  \centering
  \ximg[width=.95\linewidth,trim=10 35 30 10, clip]{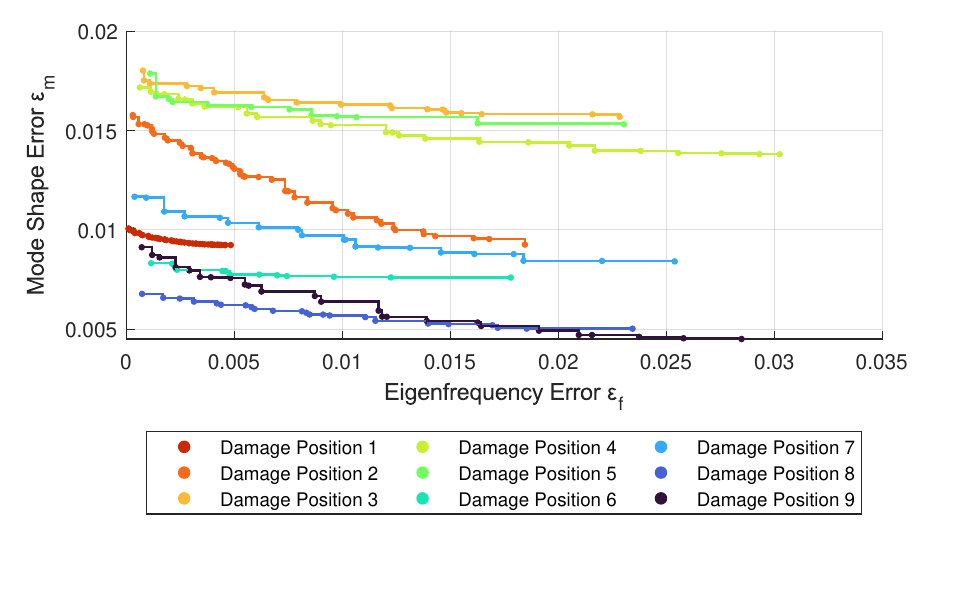}
  \caption{Non-dominated points belonging to the stiffness
    scaling factor distributions shown in \cref{fig:DamageDist}.}
  \label{fig:ParetoFrontiers}
\end{figure}

Our Algorithm \ref{alg:mogps_new} (MOGPSA) is run
with its parameters set to $T = 50$ and $N = 20$,
and with \num{1000} objective function evaluations
for each of the damage scenarios.
The resulting (first level) Pareto fronts in the objective space
are shown in \cref{fig:ParetoFrontiers}.
We observe that all nine Pareto fronts
are asymmetric with respect to the first bisector.
In particular, the fronts show small variations of the mode shape error
and large variations of the eigenfrequency error.
The number of computed non-dominated points
varies between the different damage scenarios
as it depends on the specific instance of the optimization problem.
Note that the computed points on each front
are not interpolated by a polygonal path (as usual)
but rather by concave pairs of horizontal and vertical line segments.
This guarantees that the computed fronts, including the connections,
are feasible: they belong to the extended image set $f[\Omega] + \R^2_+$,
which is bounded by the unknown true Pareto front to the lower left.

The numbers of computed points for all nine Pareto fronts
are listed in \cref{tab:cantileverXopt:results}
with minimal, average and maximal values of each damage parameter
and in \cref{tab:cantileverYopt:results}
with the corresponding values of each error measure.
Positions and values of $\mu$ and $\sigma$ refer to beam elements
and need to be multiplied by $\SI{5}{mm}$ to obtain lenght values.

\begin{table}
  \let\mc\multicolumn
  \centering
  \caption{Cantilever beam: distribution statistics of
    Pareto-optimal damage parameters for all nine scenarios.
    DP = actual damage position (beam element),
    \#\ = size of Pareto-optimal data set.}
  \label{tab:cantileverXopt:results}
  \ifcase0% Suitably rounded.
  \begin{tabular}{rr*{9}r}
   \toprule
   DP & \# & \mc3c{$\mu$} & \mc3c{$\sigma$} & \mc3c{$D$} \\
   && min & avg & max & min & avg & max & min & avg & max \\
   \midrule
    15 & 49 &  1.00&  8.40& 17.25&3.39&8.70&14.31&0.0224&0.0432&0.0563 \\
    39 & 47 & 24.75& 32.76& 37.25&4.27&8.67&15.19&0.0229&0.0332&0.0510 \\
    63 & 19 & 56.00& 58.77& 61.00&3.39&5.25& 8.16&0.0208&0.0305&0.0432 \\
    87 & 27 & 82.25& 85.08& 87.25&2.39&5.32& 8.54&0.0167&0.0281&0.0479 \\
   111 & 12 &104.75&107.25&109.75&2.13&4.52& 9.29&0.0151&0.0253&0.0391 \\
   135 & 13 &129.75&132.25&136.00&2.51&5.36& 9.16&0.0167&0.0280&0.0474 \\
   159 & 20 &156.00&157.63&159.75&2.01&4.92& 7.53&0.0146&0.0279&0.0438 \\
   183 & 23 &174.75&176.27&178.50&4.90&6.99& 9.29&0.0229&0.0327&0.0417 \\
   207 & 25 &199.75&201.40&203.50&2.76&5.54& 8.91&0.0182&0.0292&0.0495 \\
   \bottomrule
  \end{tabular}
  \or% original
  \begin{tabular}{rr*{9}r}
   \toprule
   DP & \# & \mc3c{$\mu$} & \mc3c{$\sigma$} & \mc3c{$D$} \\
   && min & avg & max & min & avg & max & min & avg & max \\
   \midrule
    15 & 49 &   1.0009 &   8.3996 &  17.253 & 3.3893 & 8.7047 & 14.309  & 0.022397 & 0.043198 & 0.05625 \\
    39 & 47 &  24.751  &  32.757  &  37.251 & 4.2679 & 8.6691 & 15.188  & 0.022918 & 0.033179 & 0.051042 \\
    63 & 19 &  56.002  &  58.765  &  61.002 & 3.3893 & 5.2456 &  8.159  & 0.020834 & 0.030483 & 0.04323 \\
    87 & 27 &  82.252  &  85.076  &  87.252 & 2.3851 & 5.3185 &  8.5355 & 0.016668 & 0.028087 & 0.047917 \\
   111 & 12 & 104.75   & 107.25   & 109.75  & 2.1341 & 4.5189 &  9.2886 & 0.015105 & 0.025261 & 0.039063 \\
   135 & 13 & 129.75   & 132.25   & 136     & 2.5106 & 5.3589 &  9.163  & 0.016667 & 0.028005 & 0.047396 \\
   159 & 20 & 156.00   & 157.63   & 159.75  & 2.0085 & 4.9205 &  7.5313 & 0.014584 & 0.027917 & 0.04375 \\
   183 & 23 & 174.75   & 176.27   & 178.5   & 4.8954 & 6.991  &  9.2885 & 0.022917 & 0.032677 & 0.041667 \\
   207 & 25 & 199.75   & 201.4    & 203.5   & 2.7615 & 5.543  &  8.9119 & 0.01823  & 0.02925  & 0.049479 \\
   \bottomrule
  \end{tabular}
  \fi
\end{table}

\begin{table}
  \let\mc\multicolumn
  \centering
  \caption{Cantilever beam: distribution statistics of
    error measures for all nine scenarios.
    DP = actual damage position (beam element).
    \#\ = size of Pareto-optimal data set.}
  \label{tab:cantileverYopt:results}
  \ifcase0% suitably rounded
  \begin{tabular}{rr*{6}r}
   \toprule
   DP & \# & \mc3c{$\varepsilon_f$} & \mc3c{$\varepsilon_m$}\\
   && min & avg & max & min & avg & max\\
   \midrule
    15 & 49 & 0.000096 & 0.0024 & 0.0048 & 0.0092 & 0.0095 & 0.01007 \\
    39 & 47 & 0.000294 & 0.0066 & 0.0185 & 0.0093 & 0.0126 & 0.01578 \\
    63 & 19 & 0.000767 & 0.0099 & 0.0228 & 0.0157 & 0.0165 & 0.01802 \\
    87 & 27 & 0.000621 & 0.0120 & 0.0303 & 0.0138 & 0.0153 & 0.01717 \\
   111 & 12 & 0.001089 & 0.0077 & 0.0230 & 0.0153 & 0.0162 & 0.01787 \\
   135 & 13 & 0.001142 & 0.0063 & 0.0178 & 0.0076 & 0.0079 & 0.00833 \\
   159 & 20 & 0.000370 & 0.0103 & 0.0254 & 0.0084 & 0.0097 & 0.01167 \\
   183 & 23 & 0.000720 & 0.0092 & 0.0234 & 0.0050 & 0.0058 & 0.00678 \\
   207 & 25 & 0.000707 & 0.0111 & 0.0285 & 0.0045 & 0.0064 & 0.00913 \\
   \bottomrule
  \end{tabular}
  \or
    \begin{tabular}{rr*{6}r}
   \toprule
   DP & \# & \mc3c{$\varepsilon_f$} & \mc3c{$\varepsilon_m$}\\
   && min & avg & max & min & avg & max\\
   \midrule
      15 & 49 & 9.6279e-05 & 0.0024017 & 0.0048286 & 0.0092413 & 0.0094939 & 0.010068 \\
      39 & 47 & 0.00029434 & 0.0065935 & 0.018451  & 0.0092648 & 0.012614  & 0.015782 \\
      63 & 19 & 0.0007667  & 0.0098939 & 0.022841  & 0.015691  & 0.016517  & 0.018017 \\
      87 & 27 & 0.00062084 & 0.012035  & 0.03025   & 0.013809  & 0.015296  & 0.017174 \\
      111 & 12 & 0.0010893  & 0.0076586 & 0.023037  & 0.015316  & 0.016165  & 0.017872 \\
      135 & 13 & 0.0011422  & 0.0062984 & 0.017796  & 0.0076057 & 0.0078725 & 0.0083276 \\
      159 & 20 & 0.00036985 & 0.010347  & 0.025388  & 0.0084186 & 0.0096959 & 0.011672 \\
      183 & 23 & 0.00072035 & 0.0092011 & 0.023434  & 0.0050329 & 0.0058271 & 0.0067849 \\
      207 & 25 & 0.00070736 & 0.011106  & 0.028483  & 0.0045128 & 0.0064417 & 0.00   \bottomrule
  \end{tabular}
  \fi
\end{table}

\begin{figure}
  \centering
  \ximg[width=.9\linewidth,trim=20 30 40 10,clip]{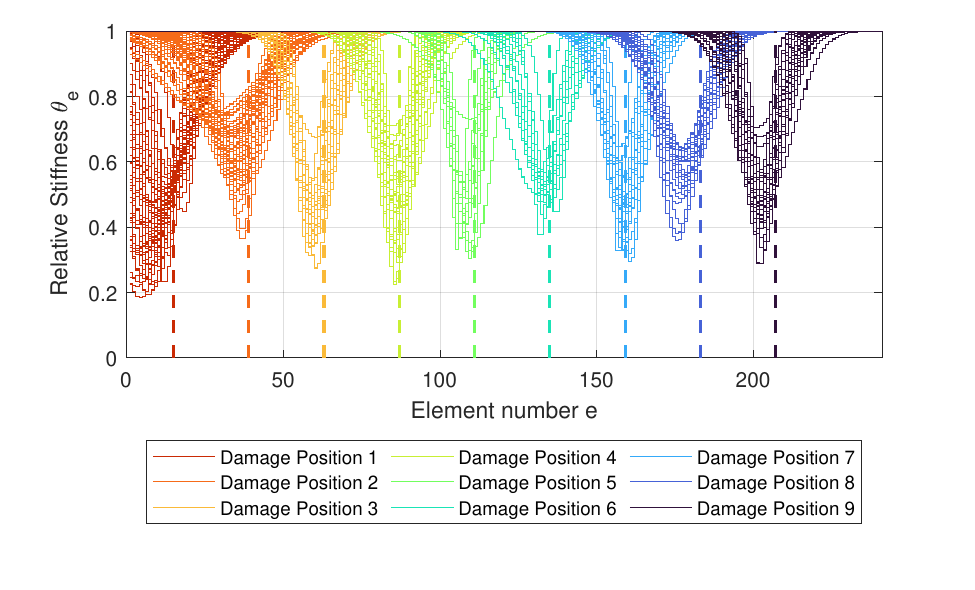}
  \caption{Plot of the stiffness scaling factor distributions identified
    for multiple different damage positions on the experimental beam.
    Scaling factors are plotted using continuous lines,
    true damage positions are plotted using dashed lines.}
  \label{fig:DamageDist}
\end{figure}

\Cref{fig:DamageDist} shows the stiffness distribution functions,
each corresponding to a point on the respective Pareto front,
for all nine damage scenarios.
The correct damage positions for all damage scenarios
are indicated with dashed lines.
The Pareto-optimal damage distributions show
that all damage positions are estimated in the correct order.
Furthermore, all damage positions
are located closely to the actual damage position
and possess a similar severity of the estimated damage.
The damage severity $D$,
which corresponds to the area enclosed by the distributions,
has similar values for most of the optimal distributions.
This is consistent with the expected result,
as all damage scenarios are implemented with the same damaged fishplate
and hence the damage severity should be approximately equal
regardless of the position.
Additionally, the widths of the depicted damage distribution functions
are consistent with each other, which also matches the actual experiment.
Minor deficiencies in the accuracy of the beam modeling
due to material irregularities and faulty model properties
can cause slight deviations in the damage position,
which can be seen throughout most given results.

Overall, these findings confirm that the model updating method
as well as the optimization algorithm
fulfill their respective tasks.
The location is adequately accurate for all scenarios considered.
All values of $D$ and $\sigma$ are well below their upper bounds
and thus well within the search area.
Additionally, the results confirm that an a-priori weighting
between eigenfrequency errors and mode shape errors
does not make sense in practice:
most of the Pareto fronts cover wide ranges of the error measures,
in particular of the eigenfrequency error,
and therefore a single point on each front
would give insufficient information.
Taking a look at Figure \ref{fig:ParetoFrontiers},
it appears that each Pareto front is
part of the boundary of a \emph{convex} extended image set $f[\Omega] + \R^2_+$,
and therefore every front also seems to be computable
by weighted sum (linear) scalarization,
\ie, by solving a series of scalarized optimization problems
which involve convex combinations of the objective functions
$\veps_f$ and $\veps_m$.
Note that for each of these problems,
the scalar objective function is also expensive
as the sum of the expensive functions $\veps_f$ and $\veps_m$.
In order to obtain a good approximation of the Pareto fronts,
the weighted sum scalarization approach would therefore require
to solve a whole series of expensive problems,
\eg, with classical pattern search algorithms.

\subsection{Girder mast structure}

\begin{figure}
  \centering
  \ximg[height=80mm]{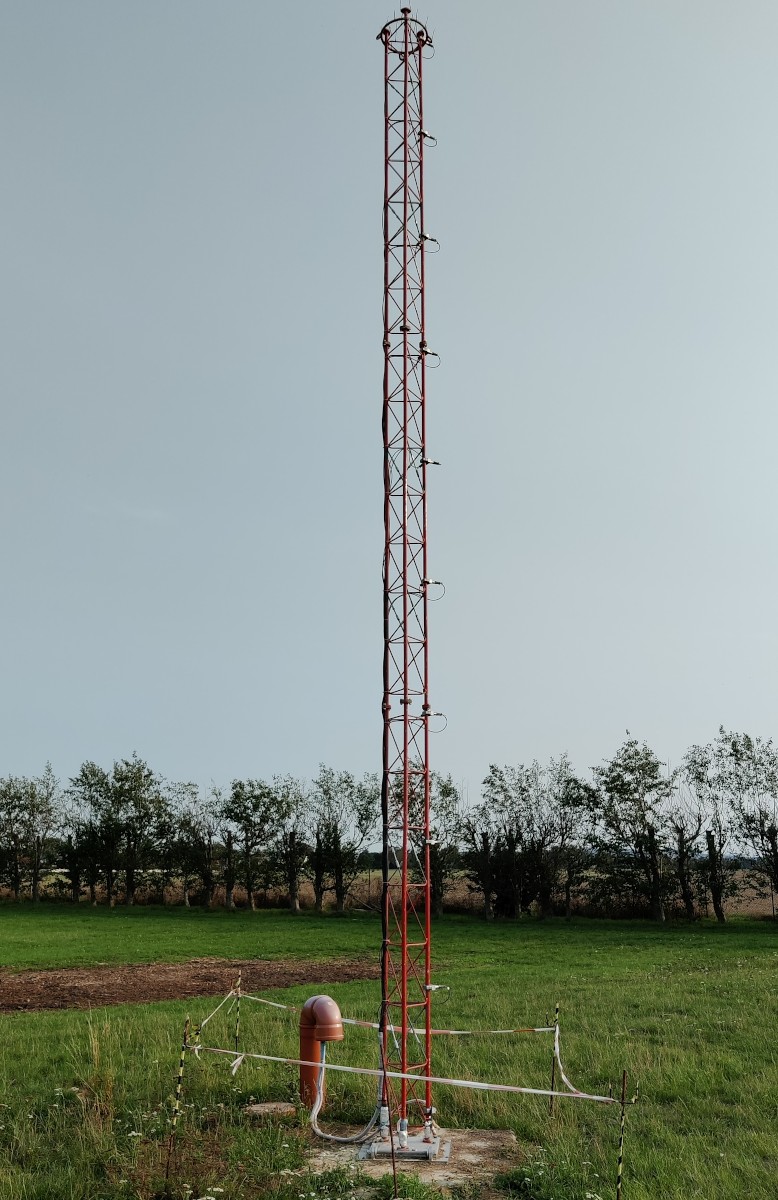}
  \hfil
  \ximg[height=80mm]{sketch_damage.eps}
  \caption{Girder mast structure: Photograph (left)
    and sketch with damage position (right).}
  \label{fig:girder}
\end{figure}

As a more complex and particularly difficult
numerical example we consider a girder mast structure.
The mast is made of construction steel
and situated in an outdoor test facility illustrated in \cref{fig:girder}.
The structure is dynamically excited by environmental conditions
such as wind forces and micro-seismic vibrations from the ground.
The whole setup is introduced in detail by Wernitz et al.\ \cite{Wernitz22}.
\begin{figure}
  \centering
  \ximg[height=60mm]{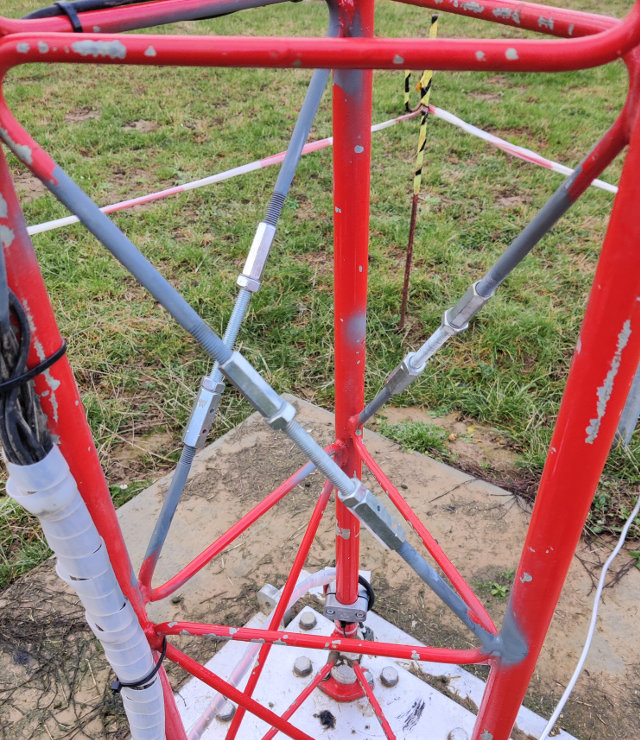}
  \caption{Photograph of reversible damage mechanisms
    integrated into the struts of the girder mast structure.}
  \label{fig:damage_strut}
\end{figure}

\subsubsection{Experimental setup}

The mast is \SI{9}{m} tall with triangular cross-section.
It consists of three vertical segments,
each made up of seven sub-segments called bays.
Each bay has three vertical cylindrical pipes
that are connected at the top and at the bottom by horizontal triangles
made of massive cylindrical cross struts with smaller diameter,
and by three inclined cross struts of the same type in between,
see \cref{fig:damage_strut}.
Each segment has three additional vertical pipes at each end,
where the segments are connected to each other without further cross struts.
The finite element model is assembled semi-automatically
in the commercial simulation software \sfname{Abaqus}.
Herein the order of beams is determined internally
in a reasonable way that keeps neigboring beams closely together.
For the girder mast there is of course no natural linear order of the beams.
All pipes (legs) and struts (braces) are modeled as Timoshenko beams
according to \cref{sec:model:struct};
they are uniformly discretized
with three elements per horizontal beam, four elements per vertical beam
and five elements per inclined beam.
Five of the six additional vertical beams that connect segments
consist of one element each, the sixth beam consists of two elements
to account for a sensor.
Altogether we have 198 beams with a total of 804 finite elements.
Finally there is a horizontal circular steel ring mounted at the mast top
(see \cref{fig:girder} left) which is modeled as a point mass.
Materials and cross-sections are assigned
according to the as-measured dimensions of the structure,
which are summarized in Table \ref{tab:girder_mast_properties}.
In this table, the nonstructural mass due to welding material,
paint and zinc coating is included in the density,
while point masses for the tower head, sensors and cables
are specified separately.
\begin{table}
  \centering
  \caption{Properties of the girder mast structure}
  \begin{tabular}{lr}
    \toprule
    Property & Value \\
    \midrule
    Elastic Modulus & \SI{210}{GPa} \\
    Outer diameter legs & \SI{30}{mm} \\
    Wall thickness legs & \SI{3}{mm} \\
    Diameter braces & \SI{10}{mm} \\
    Bay height & \SI{0.4}{m} \\
    Density including nonstructural mass & \SI{8670}{kg/m^3}\\
	Sensor node mass & \SI{0.45}{kg}\\
	Head mass & \SI{2.97}{kg}\\
	Cable mass & \SI{3.19}{kg}\\
    \bottomrule
    \end{tabular}
    \label{tab:girder_mast_properties}
\end{table}

Structural acceleration is measured
using $18$ piezoelectric accelerometers of type MEAS 8811LF-01,
situated on the nine measurement levels
indicated by ML1 to ML9 in \cref{fig:girder}.
Sensor signals are continuously converted to digital data
using National Instruments NI-9234 cards
at a sampling frequency of \SI{1651.61}{Hz},
and the data is streamed to an off-site data storage server.
In addition to the continuous measurement of the vibration,
material temperature and meteorologial data is also captured
and available for further processing.

\begin{figure}
  \centering
  \begin{subfigure}[c]{0.3\linewidth}
    \centering
    \ximg[height=58mm,trim=2 21 8 11,clip]{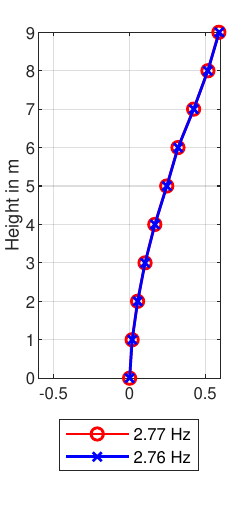}
    \subcaption{\nth{1} bending mode}
  \end{subfigure}
  \hfill
  \begin{subfigure}[c]{0.3\linewidth}
    \centering
    \ximg[height=58mm,trim=6 21 8 11,clip]{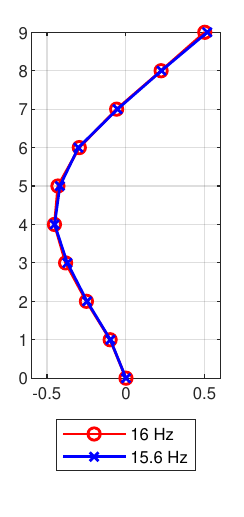}
    \subcaption{\nth{2} bending mode}
  \end{subfigure}
  \hfill
  \begin{subfigure}[c]{0.3\linewidth}
    \centering
    \ximg[height=58mm,trim=6 19 9 10,clip]{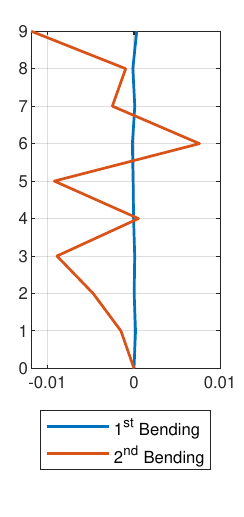}
    \subcaption{Mode difference}
  \end{subfigure}
  \caption{Mode shapes and eigenfrequencies identified from
    measurement data for the structure shown in \cref{fig:girder},
    as well as the difference between the healthy state
    and damage at position 1.
    Undamaged state indicated by red lines,
    damage at position 1 indicated by blue lines.}
  \label{fig:mode_shapes}
\end{figure}

Structural damage is introduced
by cutting braces at three distinct positions of the structure,
as indicated in \cref{fig:girder}.
The damage mechanism is depicted in \cref{fig:damage_strut}
and shows the removable struts which replace
the otherwise solid diagonal bracings.
In the damaged scenarios, a single diagonal brace is removed
at the respective position along the structure.
The modal parameters $\bPhi_k$ and $f_k$
of the structure are identified from
10-minute data sets as described in \cref{sec:model:opt}.
One of the data sets was recorded before the damage event
and the other one thereafter.
The data sets are chosen so that the vibration amplitudes
as well as the material temperature are roughly the same.
This way, data contamination due to differing
environmental conditions is prevented.
The identified modes shapes and eigenfrequencies of the structure
are shown in \cref{fig:mode_shapes}
for the undamaged state and the data set
recorded for damage position 1, respectively.
Additionally, the differences in the mode shapes with respect to
the damage are shown in a separate panel in the same figure.
\captionsetup[subfigure]{labelformat=empty}

\begin{figure}
  \centering
  \begin{subfigure}[c]{4.0em}
    \centering
    \ximg[height=65mm]{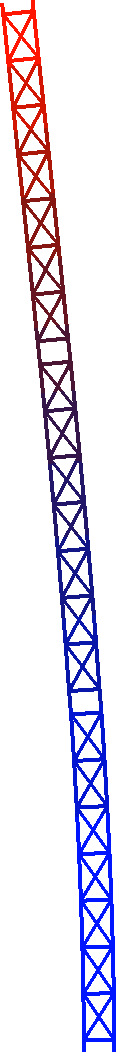}
    \subcaption{\nth{1} bending\\\SI{3.28}{Hz}}
  \end{subfigure}
  \qquad
  \begin{subfigure}[c]{4.2em}
    \centering
    \ximg[height=65mm]{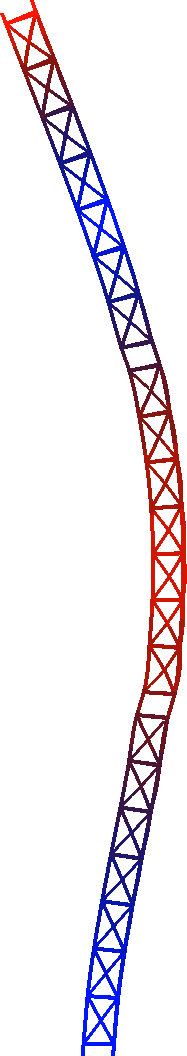}
    \subcaption{\nth{2} bending\\\SI{19.0}{Hz}}
  \end{subfigure}
  \qquad
  \begin{subfigure}[c]{7em}
    \centering
    \ximg[height=35mm]{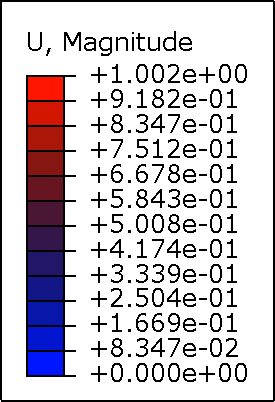}
  \end{subfigure}
  \caption{Finite element simulation results
    for mode shapes and eigenfrequencies of the undamaged state.}
  \label{fig:abq_shapes}
\end{figure}

The modal analysis of the model is carried out
using the finite element solver \sfname{Abaqus}.
The simulated mode shapes for the undamaged structure,
as shown in in \cref{fig:abq_shapes},
resemble the identified mode shapes.

% MCS: subfigure units 'em' used to match caption texts.
\captionsetup[subfigure]{labelformat=empty}

The mode shapes taken into account in the model updating procedure
belong to the first two bending mode pairs of the structure.
The eigenfrequencies drop consistently,
as expected when reducing the stiffness of the structure.
Our single-beam damage model is not entirely adequate
for the girder mast structure
because the beams are connected in a complex 3D pattern
rather than a linear sequence;
hence we need to make slight adaptations.
We assume that only a single beam is damaged
and confine the optimization to the box
$(D, \mu, \sigma) \in [0, 0.02] \x [0, 804] \x [0, 10]$.
Thus, the total length of the structure (the number of beam elements)
is our upper bound for the damage center $\mu$
whereas a fictitious longest single beam consisting of ten elements
determines the upper bound for the extent $\sigma$.
This includes a safety margin over the true maximal length of five elements.
The lower bound on the stiffness scaling factors
is enforced with parameter $\theta\tsb{min} = 0.01$.

\subsubsection{Computational results}

\begin{figure}
  \centering
  \ximg[width=0.9\linewidth,trim=10 4 37 14,clip]{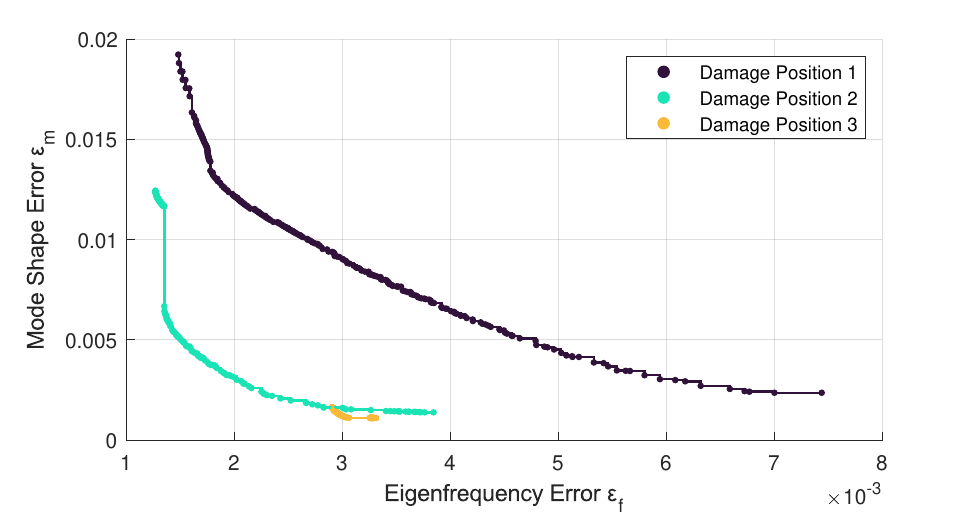}
  \caption{Pareto frontier of the girder mast model updating problem.}
  \label{fig:updating_obj}
\end{figure}

Our Algorithm \ref{alg:mogps_new} (MOGPSA) is run with
\num{5000} objective function evaluations using $T=30$ and $N=20$.
The resulting (first level) Pareto fronts
are shown in \cref{fig:updating_obj}.
The computed points on the front are again connected
by concave pairs of horizontal and vertical line segments
to visualize the space dominated by these solutions.
Overall, a large number of non-dominated points were found,
which suggests that the parameters $T$ and $N$ of the MOGPS algorithm
are adequately chosen.
The Pareto fronts have a greatly varying density of solutions
and they exhibit large gaps and concave features.
Hence, in this problem, the corresponding parts of the fronts
cannot be computed by the weighted sum (linear) scalarization technique
whereas the non-dominated sorting approach implemented in our algorithm
provides all the information required in practical applications.

The minimal, median and maximal values of each damage parameter
and each error measure, respectively,
are listed in \cref{tab:mastxopt:results,tab:mastyopt:results}.
Positions and values of $\mu$ and $\sigma$ refer again to beam elements.
As mentioned, the order of
beams is determined by the \sfname{Abaqus} software,
and here it is arbitrary since the beams of the mast
do not form a linear sequence.
We choose the medians of the distributions here rather than the average
since the computed Pareto optimal points form two clusters
(in contrast to a single cluster for the steel beam).
These clusters turn out to correspond to respective subsets of points
associated with correct and incorrect damage locations, see next paragraph.

\begin{table}
  \let\mc\multicolumn
  \centering
  \caption{Girder mast: distribution statistics of
	    Pareto-optimal damage parameters for all three scenarios.
	    DP = actual damage position (beam element),
	    \#\ = size of Pareto-optimal data set.}
  \label{tab:mastxopt:results}
  \begin{tabular}{rr*{9}r}
	   \toprule
	   DP & \# & \mc3c{$\mu$} & \mc3c{$\sigma$} & \mc3c{$D \x 1000$} \\
	   && min & median & max & min & median & max & min & median & max \\
	   \midrule
	      6 & 262 &   3.48&  9.63&163.06&4.68&6.34&7.88&14.6&15.9& 18.3\\
	    274 & 179 & 273.12&276.51&358.72&1.95&2.27&4.73& 6.1& 6.8& 13.7\\
	    429 & 667 &   1.01&437.27&   804&1e-4&0.11&  10&   0&0.24&  4.5\\
	   \bottomrule
	\end{tabular}
\end{table}

\begin{table}
	\let\mc\multicolumn
	\centering
	\caption{Girder mast: distribution statistics of
		error measures for all three scenarios.
		DP = actual damage position (beam element),
		\#\ = size of Pareto-optimal data set.}
	\label{tab:mastyopt:results}
	\begin{tabular}{rr*{6}r}
		\toprule
		DP & \# & \mc3c{$\varepsilon_f$} & \mc3c{$\varepsilon_m$}\\
		&& min & median & max & min & median & max\\
		\midrule
		6 & 262 & 0.0015 & 0.0026 & 0.0074 & 0.0024 & 0.0103 & 0.0192 \\
		274 & 179 & 0.0013 & 0.0014 & 0.0038 & 0.0014 & 0.0057 & 0.0124 \\
		429 & 667 & 0.0029 & 0.0029 & 0.0033 & 0.0011 & 0.0016 & 0.0016 \\
		\bottomrule
	\end{tabular}
\end{table}

\begin{figure}
  \centering
  \ximg[width=0.9\linewidth,trim=20 3 38 15,clip]{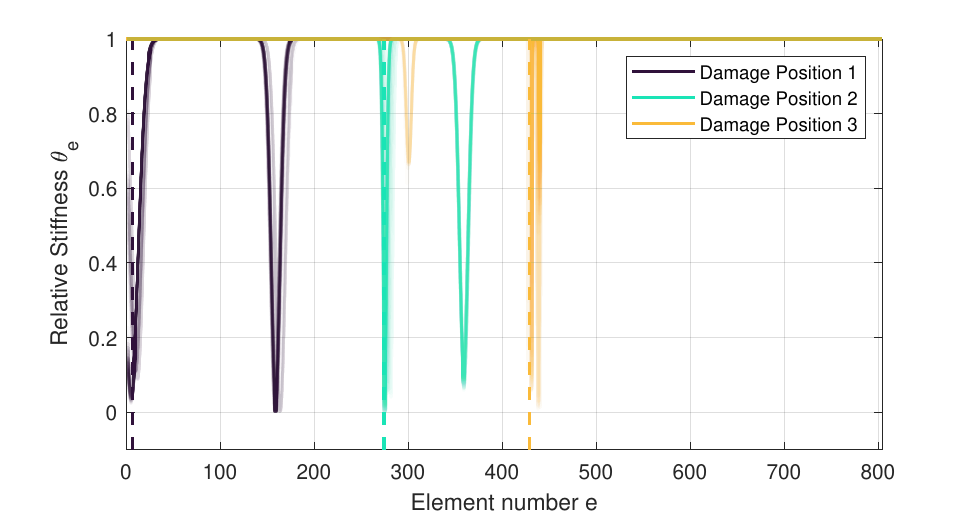}
  \caption{Stiffness scaling factors computed for the Pareto fronts.
    MOGPS results are shown using solid lines and position
    of experimentally damaged bay is indicated using dashed lines.}
  \label{fig:updating_mode_shapes}
\end{figure}

The non-dominated stiffness distributions
for all three damage positions
are shown in \cref{fig:updating_mode_shapes}.
Dashed lines in the corresponding colors indicate the exact damage positions.
Each colored curve shows an overlay of several curves
corresponding to the individual non-dominated points.
A major part of the optimal damage distributions
show the estimated damage position close to the actual damage position:
121 for damage position one, 83 for damage position two,
and 573 for damage position three.
However, we also observe 141, 96 respectively 94
optimal damage distributions that identify the wrong beam as damaged.

The obtained results show
that the presented model updating procedure
enables a precise damage location
for all investigated damage positions.
A fraction of solutions for all damage positions
seem to identify the damage position far from the actual position.
This is due to the bracings in the finite elements
not being sorted in any particular order
when applying the damage distribution.
In fact, the addressed damage distributions all identify the maximum damage
on diagonal braces close to the correct position,
since neighboring damage positions tend to have a similar impact
on the dynamic behavior of the structure.
Finally, most of the Pareto-optimal solutions
posses a minimum relative stiffness close to zero,
which is consistent with the experiment,
where the diagonal braces are completely removed.
This further strengthens the results
as all of the computed damage distributions
have minimum relative stiffnesses near zero.

%%% Local Variables:
%%% mode: latex
%%% TeX-master: t
%%% End:

\section{Conclusion}

\label{sec:conclusion}
The application of the multi-objective model updating scheme
to the two experimental structures yields satisfactory results overall.
The proposed pattern-search algorithm is capable
of finding non-dominated solutions effectively
and structural damage is located adequately.
The results for the laboratory beam experiment are almost perfect.
The damage location in the girder mast example is ambiguous
since a fraction of the non-dominated solutions
point to damage positions next to the actual damage.
This is not surprising as the one-dimensional damage model
cannot adequately capture the complex three-dimensional
connection topology of the individual beams.
In spite of this, our approach does detect the presence of a damage reliably
and narrows down the damage center to just two possible locations,
which is still good enough from a practical point of view.

Regarding the proposed optimization algorithm,
there is room for future research.
A subsequent theoretical analysis should address questions about
the approximation quality of computed Pareto optimal points
and the quality of coverage of the Pareto front.
Further, numerical benchmarks using academic test problems
and additional application problems should be conducted
to evaluate the performance of the algorithm,
to study the influence of the parameters $T$ and $N$,
and to compare it to other algorithms from the literature.

%%% Local Variables:
%%% mode: latex
%%% TeX-master: t
%%% End:

\section*{Acknowledgement}

The third author wishes to thank for
the financial support by the Deutsche Forschungsgemeinschaft
(DFG, German Research Foundation) -- SFB1463 -- 434502799, subproject C04.
The fifth author is funded by the
Federal Ministry for Economic Affairs and Energy
\emph{MMRB-Repair-care -- Multivariate damage monitoring of rotor blades:
  implementation and analysis of the effects of repair measures},
FKZ 03EE2043C.
The sixth and seventh authors gratefully acknowledge
the financial support by the Deutsche Forschungsgemeinschaft
-- SFB1463 -- 434502799, subproject B05.

\bibliographystyle{siam}
\bibliography{literature}

\end{document}